\documentclass[12pt]{article}

\usepackage{graphicx}
\usepackage{caption}
\usepackage[margin=1in]{geometry}

\usepackage{amstext}
\usepackage{latexsym}
\usepackage{amssymb}
\usepackage{geometry}

\def\R{\mathbb R}
\def\C{\mathbb C}
\def\N{\mathbb N}

\newcommand{\eps}{\varepsilon}

\begin{document}

\newtheorem{theorem}{Theorem}[section]
\renewcommand{\thetheorem}{\arabic{section}.\arabic{theorem}}
\newtheorem{definition}[theorem]{Definition}
\newtheorem{deflem}[theorem]{Definition and Lemma}
\newtheorem{lemma}[theorem]{Lemma}
\newtheorem{example}[theorem]{Example}
\newtheorem{remark}[theorem]{Remark}
\newtheorem{remarks}[theorem]{Remarks}
\newtheorem{cor}[theorem]{Corollary}
\newtheorem{pro}[theorem]{Proposition}
\newtheorem{proposition}[theorem]{Proposition}

\renewcommand{\theequation}{\thesection.\arabic{equation}}

\title{Bifurcation from the Kurth solution in galactic dynamics}
\author{{\sc Markus Kunze$^{1}$ \& Rafael Ortega$^{2}$} \\[2ex]
       $^{1}$ Universit\"at K\"oln, Institut f\"ur Mathematik, Weyertal 86-90, \\
       D\,-\,50931 K\"oln, Germany \\[1ex]
       $^{2}$ Departamento de Matem\'{a}tica Aplicada, Universidad de Granada, \\
       E-18071 Granada, Spain \\[1ex]
       }
\date{}
\maketitle
\begin{abstract}\noindent 
It will be shown that there exists an infinite-dimensional continuum
${\cal C}$ of weak static solutions of the Vlasov-Poisson system 
that bifurcates from the Kurth solution. Each $f_\ast\in {\cal C}$ 
has the charge density $\rho_{f_\ast}=\rho_{{\rm Kurth}}$, 
and (like the Kurth solution itself) each $f_\ast$ is surrounded 
by time-periodic weak solutions. 
\end{abstract}

%%%%%%%%%%%%%%%%%%%%%%%%%%%%%%%%%%%%%%%%%%%%%%%%%%%%%%%%%%%%%%%%%%%%%%%%%%%%%%%%%%%%

\setcounter{equation}{0}
\section{Introduction and main results}

The gravitational Vlasov-Poisson system is a well-established model 
for the evolution of a self-gravitating, non-relativistic ensemble 
of a large number of stars, galaxies or even clusters of galaxies \cite{BT}. 
From a mathematical viewpoint, around 1990 the existence and uniqueness theory 
of solutions has been completed in a very satisfactory manner \cite{pfaff,LP,schaeff}. 
In addition, a great deal is known about the existence and stability of steady states, 
see \cite{LMR12,GRking} and the references therein. However, 
apart from that, there is close to no rigorous result 
for the unmodified and fully nonlinear system on $\R^3\times\R^3$ 
near isotropic or anisotropic steady states (like the Kurth solution or the polytropes) 
that has a dynamical systems flavor. For instance, 
it would be very interesting to address the local behavior close to (an-)isotropic polytropes; 
see \cite{RaRe,Str} for some numerical results, which seem hard to obtain analytically. 
\smallskip 

\noindent 
On a linear level, properties of the linearization about such (an-)isotropic equilibrium 
are now better understood \cite{K_bk, HRS, Wed}, and in some works a point mass at the center is added 
and perturbations of the linear system off this point mass 
are considered \cite{damp_vers_osc, HS}. 
\smallskip 

\noindent 
For the nonlinear system, so far mostly the long-time behavior and/or (modified) scattering 
has been treated for small solutions \cite{ChK, IPWW2, FOPW}, close to certain homogeneous equilibria 
(depending on $|v|$ only) \cite{IPWW} in the plasma case,  
and under the addition of a Kepler potential \cite{CL} 
or an attracting point mass \cite{KeWi}. 
\medskip

In the present paper, our starting point will be the Kurth solution 
$f_{{\rm Kurth}}=f_{{\rm Kurth}}(x, v)$ from \cite{RK}, also see \cite{RaRe,HRS,2ndlook} 
for more context. It is a particular (weak) static solution, 
which is surrounded by time-periodic solutions that are almost explicit.  
A particular feature of this solution is that it gives rise 
to a constant density $\rho_{{\rm Kurth}}=\rho_{{\rm Kurth}}(x)$, 
which is a multiple of the characteristic function of the unit ball in $\R^3$. 
We will construct an infinite-dimensional continuum
${\cal C}$ of weak static solutions 
that bifurcates from $f_{{\rm Kurth}}$ in a well-defined sense. 
Each $f_\ast\in {\cal C}$ has the charge density $\rho_{f_\ast}=\rho_{{\rm Kurth}}$, 
and each $f_\ast$ is surrounded by time-periodic weak solutions as well. 
Hence we are able to obtain an abundance 
of new solutions, together with a mind-boggling bifurcation scenario. 
Both the static solutions and the time-periodic solutions are spherically symmetric. 

Our method of proof is quite explicit, and it relies on the uniform distribution 
of the bifurcating solutions. Therefore it is not immediately clear 
if it can be adapted to e.g.~the polytropes. However, there are other classes 
of steady states, for which the strategy could be attempted. In fact the Kurth solution 
is the first member of a family of equilibria described in \cite{Ahmad}, 
which are also called Osipkov-Merritt models \cite{BT}. 
For instance, its second member has the density $\rho(x)=(1-|x|^2)_+$, 
corresponding to a potential of order $4$ for $|x|\le 1$, 
and it could be a natural candidate to try to prove bifurcation from.   
\medskip 

To add some more details on the Kurth solution, we recall the Vlasov-Poisson system 
in the gravitational case as
\begin{equation}\label{vpgr1} 
   \partial_t f(t, x, v)+v\cdot\nabla_x f(t, x, v)
   -\nabla_x U_f(t, x)\cdot\nabla_v f(t, x, v)=0,
\end{equation} 
where 
\begin{equation}\label{vpgr2} 
   \Delta U_f(t, x)=4\pi\rho_f(t, x),\quad\lim_{|x|\to\infty} U_f(t, x)=0,
   \quad\rho_f(t, x)=\int_{\R^3} f(t, x, v)\,dv,
\end{equation} 
for $(t, x, v)\in\R\times\R^3\times\R^3$. Therefore 
\begin{equation}\label{vpgr2b} 
   U_f(t, x)=-\int_{\R^3}\frac{\rho_f(t, y)}{|x-y|}\,dy.
\end{equation} 
Note that the only unknown is the scalar particle density function $f=f(t, x, v)\ge 0$. 

The Kurth solution is a weak static solution and given by 
\begin{eqnarray}\label{fKdef} 
   & & f_{{\rm Kurth}}(x, v)=\frac{3}{4\pi^3}\,
   \frac{1}{(1-|x|^2-|v|^2+|x\wedge v|^2)^{1/2}}
   \quad\mbox{where}\quad (\ldots)>0
   \,\,\mbox{and}\,\,|x\wedge v|<1,
   \nonumber
   \\[1ex] & & \mbox{and}\quad f_{{\rm Kurth}}(x, v)=0\quad\mbox{else}, 
\end{eqnarray} 
for $(x, v)\in\R^3\times\R^3$. 

\begin{remark}  
{\rm Note that the inequality $1-|x|^2-|v|^2+|x\wedge v|^2>0$ does not imply 
that $|x\wedge v|<1$. This will be very clear after the geometric 
discussions of Section \ref{supp_sect}. Also, it is easy to find explicit 
$(x, v)$ satisfying the first inequality, but not the second.} 
\end{remark} 

The charge density associated to (\ref{fKdef}),  
\begin{equation}\label{hano} 
   \rho_{{\rm Kurth}}(x)=\int_{\R^3} f_{{\rm Kurth}}(x, v)\,dv
   =\frac{3}{4\pi}\,{\bf 1}_{\overline{B_1(0)}}\,(x), 
\end{equation}  
is, up to a normalizing factor, the characteristic function of the unit ball in $\R^3$. 
The solution to $\Delta U_{{\rm Kurth}}=4\pi\rho_{{\rm Kurth}}$ 
and $U_{{\rm Kurth}}(x)\to 0$ as $|x|\to\infty$ is given by 
\begin{equation}\label{UQkurth} 
   U_{{\rm Kurth}}(x)=\left\{\begin{array}{c@{\quad:\quad}c} 
   \frac{1}{2}\,|x|^2-\frac{3}{2} & |x|\le 1
   \\[1ex] -\frac{1}{|x|} & |x|>1
   \end{array}\right. ;
\end{equation}
also cf.~\cite{2ndlook} for further details. 

A remarkable special feature of the Kurth solution 
is that it is surrounded by time-periodic solutions 
$f_{{\rm Kurth}, \eps}(t, x, v)$, as discovered in \cite{RK}, 
whose period $T_\eps\to 2\pi$ as $\eps\to 0$. To introduce those solutions, 
consider the Hamiltonian system $H(\phi, \dot{\phi})
=\frac{1}{2}\,\dot{\phi}^2+V(\phi)$ for $V(\phi)
=-\frac{1}{\phi}+\frac{1}{2\phi^2}$. The associated equation of motion 
is $\ddot{\phi}=-\frac{1}{\phi^2}+\frac{1}{\phi^3}$. Let $\phi_\eps$ 
denote the solution such that 
\begin{equation}\label{phieps_IVP} 
   \phi_\eps(0)=1\quad\mbox{and}\quad\dot{\phi}_\eps(0)=\eps.
\end{equation}   
Then $\phi_\eps$ is periodic for $\eps<1$ 
(in fact $|\eps|<1$), and the period is calculated to be
\begin{equation}\label{kurth_per} 
   T_\eps=2\int_{\phi_{{\rm min}}^\eps}^{\phi_{{\rm max}}^\eps}
   \frac{ds}{\sqrt{2(E_\eps-V(s))}}
   =\frac{2\pi}{(1-\eps^2)^{3/2}}
\end{equation}  
for $E_\eps=H(\phi_\eps, \dot{\phi}_\eps)=-\frac{1}{2}(1-\eps^2)
\in [-\frac{1}{2}, 0[$, $\phi_{{\rm min}}^\eps=\frac{1}{1+\eps}$, 
and $\phi_{{\rm max}}^\eps=\frac{1}{1-\epsilon}$. 
In particular, $T_\eps\to 2\pi$ as $\eps\to 0$. Then 
\[ f_{{\rm Kurth}, \eps}(t, x, v)=f_{{\rm Kurth}}
   \Big(\frac{x}{\phi_\eps(t)}, \phi_\eps(t)v-\dot{\phi}_\eps(t)x\Big),
   \quad t\in\R,\quad (x, v)\in\R^3\times\R^3. \]   
The associated density is 
\begin{equation}\label{rhoeps} 
   \rho_{{\rm Kurth}, \eps}(t, x)
   =\int_{\R^3} f_{{\rm Kurth}, \eps}(t, x, v)\,dv
   =\frac{3}{4\pi}\,\frac{1}{\phi_\eps(t)^3}\,{\bf 1}_{\{|x|<\phi_\eps(t)\}}
   =\frac{1}{\phi_\eps(t)^3}\,\rho_{{\rm Kurth}}
   \Big(\frac{x}{\phi_\eps(t)}\Big),
\end{equation}  
resulting in the potential 
\begin{equation}\label{Ueps} 
   U_{{\rm Kurth}, \eps}(t, x)
   =\frac{1}{\phi_\eps(t)}\,U_{{\rm Kurth}}\Big(\frac{x}{\phi_\eps(t)}\Big).
\end{equation} 

Our main results are as follows. 
See Definition \ref{weak_sol_def} below for the notion of a weak solution 
that we are going to use. 

\begin{theorem}\label{ex_fGam} 
For every $\Gamma\in [1-2^{-12}, 1[$, there exists 
a weak static solution $f_\Gamma\in L^1(\R^3\times\R^3)$ 
whose essential support is 
\begin{eqnarray*} 
   {\cal D}_\Gamma & = & \{(x, v)\in\R^3\times\R^3: 
   -1\le |x\wedge v|^2-|x|^2-|v|^2\le 0, |x\wedge v|^2\le\Gamma\}
   \\ & \subset & \{(x, v)\in\R^3\times\R^3:  |x|\le 1, |v|\le 1\} 
\end{eqnarray*}  
and which satisfies 
\[ \rho_{f_\Gamma}(x)=\frac{3}{4\pi}\,{\bf 1}_{\overline{B_1(0)}}\,(x)
   =\rho_{{\rm Kurth}}(x)
   \quad\mbox{for a.e.}\quad x\in\R^3. \]
\end{theorem} 

The proof of Theorem \ref{ex_fGam} is given in Section \ref{PhiGam_sect}, 
along with the proof of the following corollary, 
which gives more information about the discontinuity set of $f_\Gamma$. 

\begin{cor}\label{sing_cor} Let 
\begin{equation}\label{calM_def} 
   {\cal M}_\Gamma
   =\{(x, v)\in\R^3\times\R^3: |x|^2+|v|^2-|x\wedge v|^2=1, 
   |x\wedge v|^2\le\Gamma\}\subset {\cal D}_\Gamma
\end{equation}   
and 
\begin{equation}\label{calN_def} 
   {\cal N}_\Gamma
   =\{(x, v)\in\R^3\times\R^3: |x|^2+|v|^2-|x\wedge v|^2=\Gamma, 
   |x\wedge v|^2\le\Gamma\}\subset {\cal D}_\Gamma.
\end{equation}  
Then $f_\Gamma$ is continuous and positive on 
${\cal D}_\Gamma\setminus ({\cal M}_\Gamma\cup {\cal N}_\Gamma)$. 
Moreover, there exists $l_\Gamma>0$ such that 
\begin{equation}\label{cor_c1} 
   f_\Gamma(x, v)\to l_\Gamma\quad\mbox{as}
   \quad (x, v)\searrow {\cal N}_\Gamma
\end{equation}  
and furthermore 
\begin{equation}\label{cor_c2} 
   f_\Gamma(x, v)\to\infty\quad\mbox{as}
   \quad (x, v)\nearrow{\cal N}_\Gamma
   \,\,\mbox{or}\,\,(x, v)\nearrow {\cal M}_\Gamma.
\end{equation}  
\end{cor} 
The notation $(x, v)\searrow {\cal N}_\Gamma$ indicates 
convergence to ${\cal N}_\Gamma$ from ``above'', i.e., 
\[ (x, v)\in {\cal D}_\Gamma,
   \quad |x|^2+|v|^2-|x\wedge v|^2>\Gamma,
   \quad {\rm dist}((x, v), {\cal N}_\Gamma)\to 0. \] 
The other cases are defined similarly. 
From our proof we will see that $l_\Gamma\to\infty$ 
as $\Gamma\nearrow 1$. 
\medskip 

Next we will present a result 
that makes precise in which sense bifurcation 
from the Kurth solution does occur. 

\begin{theorem}\label{Gamto1} For $\Gamma\in [1-2^{-12}, 1[$ 
let $f_\Gamma$ be as in Theorem \ref{ex_fGam}. Then 
\[ f_\Gamma\to f_{{\rm Kurth}}\quad\mbox{in}\quad L^1(\R^3\times\R^3)  
   \quad\mbox{as}\quad\Gamma\nearrow 1. \]  
\end{theorem} 

We are also able to show that, 
like in the case of the Kurth solution, the new weak static solutions 
$f_\Gamma$ are surrounded by periodic solutions. 

\begin{theorem}\label{per_gam} For $\Gamma\in [1-2^{-12}, 1[$ 
let $f_\Gamma$ be as in Theorem \ref{ex_fGam}. Let $\phi=\phi(t)$ 
for $t\in I$ be a solution to 
\begin{equation}\label{Kep}
   \ddot{\phi}=-\frac{1}{\phi^2}+\frac{1}{\phi^3},
   \quad\phi>0.  
\end{equation}
Then 
\[ f_{\Gamma, \phi}(t, x, v)
   =f_\Gamma
   \Big(\frac{x}{\phi(t)}, \phi(t)v-\dot{\phi}(t)x\Big),
   \quad t\in I,\quad (x, v)\in\R^3\times\R^3, \] 
is a weak solution. In particular, if $\phi_\eps$ denotes the solution 
to (\ref{Kep}) for $\eps<1$ such that (\ref{phieps_IVP}) holds, then 
\[ f_{\Gamma, \eps}(t, x, v)
   =f_\Gamma
   \Big(\frac{x}{\phi_\eps(t)}, \phi_\eps(t)v-\dot{\phi}_\eps(t)x\Big),
   \quad t\in\R,\quad (x, v)\in\R^3\times\R^3, \] 
is a weak solution that is $T_\eps$-periodic in $t$, 
where $T_\eps$ is from (\ref{kurth_per}).  
\end{theorem}  

The proof of Theorem \ref{per_gam} is elaborated in Section \ref{Liling}, 
where we also have a remark (Remark \ref{Kep_rem}) that explains 
in more detail the origin of (\ref{Kep}). 
\medskip

Finally, in Section \ref{infinite_sect} we can combine the previous results 
to obtain a plethora of other static and time-periodic weak solutions 
by convexity arguments. 

\begin{theorem}\label{infi_thm} There exists an infinite-dimensional continuum
${\cal C}$ of weak static solutions. Each $f_\ast\in {\cal C}$ 
has the charge density $\rho_{f_\ast}=\rho_{{\rm Kurth}}$, 
and $f_\ast$ is surrounded by time-periodic weak solutions. 
\end{theorem} 
\medskip

We are going to close this introduction with the definition of a weak solution 
that we will be using throughout. 

\begin{definition}[weak solution]\label{weak_sol_def} 
In $L^1(\R^3\times\R^3)$ denote by 
\[ L^1_+(\R^3\times\R^3)=\{f\in L^1(\R^3\times\R^3): f\ge 0\,\,{\rm a.e.}\} \]  
the cone of a.e.~non-negative functions $f=f(x, v)$. 
We will consider a fixed open time interval $I\subset\R$; 
arbitrary compact subintervals will be written as 
$J\subset\joinrel\subset I$. 
\smallskip 

\noindent
A continuous function $f: I\to L^1_+(\R^3\times\R^3)$ 
will be called a weak solution, if 
\begin{itemize} 
\item[(a)] for every $J\subset\joinrel\subset I$ there exists a compact set 
$K_J\subset\R^3\times\R^3$ such that, for each $t\in J$, 
$f(t, x, v)=0$ for a.e.~$(x, v)\in (\R^3\times\R^3)\setminus K_J$; 
\item[(b)] $\rho_f\in L^\infty(J\times\R^3)$ 
for every $J\subset\joinrel\subset I$, where 
\[ \rho_f(t, x)=\int_{\R^3} f(t, x, v)\,dv,
   \quad t\in I,\quad {\rm a.e.}\,\,x\in\R^3; \] 
\item[(c)] $f$ is a distributional solution 
of the Vlasov equation (\ref{vpgr1}), i.e., it holds that 
\begin{eqnarray}\label{distr_eq} 
   & & \int_I\int_{\R^3}\int_{\R^3} f(t, x, v)
   \,\Big(\partial_t\varphi(t, x, v)+v\cdot\nabla_x\varphi(t, x, v)
   -\nabla_x U_f(t, x)\cdot\nabla_v\varphi(t, x, v)\Big)\,dt\,dx\,dv
   \nonumber   
   \\ & & \,\,=\,\,0 
\end{eqnarray} 
for all functions $\varphi\in {\cal D}(I\times\R^3\times\R^3)$, where 
\[ U_f(t, x)=-\int_{\R^3}\frac{\rho_f(t, y)}{|x-y|}\,dy,
   \quad t\in I,\quad {\rm a.e.}\,\,x\in\R^3. \]  
\end{itemize}
If $f$ is independent of time, it is called a weak static solution. 
In that case we let $I=\R$, and test functions can be taken in the class 
${\cal D}(\R^3\times\R^3)$. 
\end{definition} 

\begin{remark}  
{\rm This definition is more or less standard 
and adapted from, for instance, \cite[Definition 1.1]{loeper}. 
Since we will mostly be interested in periodic solutions, 
the initial data at $t=0$ do not play a role. 
Concerning the integrability of the nonlinear term in (\ref{distr_eq}), 
we recall from \cite[Thm.~10.2(iii)]{LiLo} the following facts 
about Poisson's equation: Given $\rho\in L^\infty(\R^3)$ 
with compact support, the potential 
\[ U(x)=-\int_{\R^3}\frac{\rho(y)}{|x-y|}\,dy,
   \quad x\in\R^3, \] 
is in $C^1(\R^3)$ and satisfies  
\[ \nabla U(x)=\int_{\R^3}\frac{\rho(y)}{|x-y|^3}\,(x-y)\,dy,
   \quad x\in\R^3. \]
In particular, it follows that 
\[ {\|\nabla U\|}_{L^\infty(\R^3)}
   \le C {\|\rho\|}_{L^\infty(\R^3)}, \]  
where $C>0$ only depends on the size of the support of $\rho$. 
Coming back to (\ref{distr_eq}), in view of the assumptions 
(a) and (b) this observation can be applied to $U_f(t, \cdot)$. 
Then $\nabla U_f$ is a measurable function 
and $\nabla_x U_f\in L^\infty(J\times\R^3)$ 
for every $J\subset\joinrel\subset I$. 
As a consequence, we obtain $f\nabla_x U_f\cdot\nabla_v\varphi
\in L^1(I\times\R^3\times\R^3)$ for every test function 
$\varphi\in {\cal D}(I\times\R^3\times\R^3)$. 
}  
\end{remark} 

%%%%%%%%%%%%%%%%%%%%%%%%%%%%%%%%%%%%%%%%%%%%%%%%%%%%%%%%%%%%%%%%%%%%%%%%%%%%%%%%%%

\setcounter{equation}{0}
\section{The support of the steady states}
\label{supp_sect} 

Before actually going into the details about the new steady states 
in Section \ref{constr_sect}, we need to include a few considerations 
concerning their support. Suppose that $Q=Q(u): \R\to [0, \infty[$ 
is an integrable function with support in $[-1, 0]$. 
We would like to determine $Q$ in such a way that 
\begin{equation}\label{fGam} 
   f(x, v)=Q(\ell^2(x, v)-|x|^2-|v|^2)\,{\bf 1}_{\{\ell^2(x, v)\le\Gamma\}}
\end{equation}  
is a weak static solution and has $\rho_f=\rho_{{\rm Kurth}}$; 
here $\ell(x, v)=|x\wedge v|$. Observe that 
\[ \ell^2(x, v)-|x|^2-|v|^2=\ell^2(x, v)\Big(1-\frac{1}{r^2}\Big)-r^2-p_r^2 \] 
in the coordinates $(r, p_r, \ell)\in [0, \infty[\times\R\times [0, \infty[$ 
for $r=|x|$ and $p_r=\frac{x\cdot v}{r}$, cf.~the appendix, Section \ref{append_sect}. Then
\begin{equation}\label{picc} 
   -1\le |x\wedge v|^2-|x|^2-|v|^2\le 0,\quad |x\wedge v|^2\le\Gamma, 
\end{equation}  
is rewritten as 
\begin{equation}\label{pic2c} 
   0\le r^2+p_r^2+\frac{\ell^2}{r^2}(1-r^2)\le 1,\quad\ell^2\le\Gamma. 
\end{equation}

The first observation to make is that (\ref{picc}) or (\ref{pic2c}) implies $|x|\le 1$. 
To see this, we note the relation 
\begin{equation}\label{mahe} 
   |x|^2+|v|^2\le 1+|x\wedge v|^2\le 1+|x|^2|v|^2,
\end{equation} 
and therefore $(1-|x|^2)(1-|v|^2)\ge 0$. Thus if we assume $|x|=1+\delta$ 
for some $\delta>0$, then we must have $|v|\ge 1$, But then 
\[ (1+\delta)^2+1\le |x|^2+|v|^2\le 1+|x\wedge v|^2\le 1+\Gamma\le 2, \] 
which is a contradiction. Hence we need to consider only $r\in [0, 1]$ 
in what follows, and the requirement $\ldots\ge 0$ could be dropped from (\ref{pic2c}).    
Additionally, due to (\ref{mahe}) the conditions (\ref{picc}) or (\ref{pic2c}) 
also imply $|v|\le 1$. 
\smallskip

Later it will be more convenient to use the variables $(r, p_r, u)$, where 
\begin{equation}\label{ellu} 
   -u=r^2+p_r^2+\frac{\ell^2}{r^2}(1-r^2).
\end{equation}  
Then (\ref{pic2c}) reads as 
\begin{equation}\label{pic2d} 
   -1\le u\le 0,\quad 0\le -r^2(u+r^2+p_r^2)\le (1-r^2)\Gamma. 
\end{equation} 

The following support restrictions arise from (\ref{pic2c}) or (\ref{pic2d}). 

\begin{lemma}\label{Gamsuppc} Let 
\[ D(r)=\{(p_r, u)\in\R\times\R: 
   (\ref{pic2c})\,\,{\rm holds\,\,for}\,\,(r, p_r, \ell^2)\}. \] 
\begin{itemize} 
\item[(a)] If $r>1$, then $D(r)=\emptyset$. 
\item[(b)] If $0\le r^2\le\Gamma$, then 
\begin{equation}\label{albr1} 
   D(r)=\{(p_r, u)\in\R\times [-1, 0]: 
   -1\le u\le -r^2-p_r^2\}.
\end{equation}  
\item[(c)] If $\Gamma\le r^2\le 1$, then 
\begin{equation}\label{albr2}  
   D(r)=\{(p_r, u)\in\R\times [-1, 0]: 
   -\gamma-p_r^2\le u\le -r^2-p_r^2\},
\end{equation}  
where for $r\in ]0, 1]$ and $\Gamma\in ]0, 1[$ we abbreviate 
\begin{equation}\label{gamma_def} 
   \gamma=r^2+\Big(\frac{1-r^2}{r^2}\Big)\,\Gamma.
\end{equation} 
\end{itemize} 
\end{lemma} 
{\bf Proof\,:} (a) has been shown above. 
To establish (b) and (c), we assume $r\in [0, 1]$ 
and observe that an obvious algebraic manipulation shows that 
(\ref{pic2d}) is equivalent to 
\[ u\ge -1,\quad -\gamma-p_r^2\le u\le -r^2-p_r^2. \] 
Since $\gamma<1$ if and only if $r^2\in ]\Gamma, 1[$, 
the conclusion follows. 
{\hfill$\Box$}\bigskip 

\noindent
{\bf Remark:} In the $(p_r, u)$-plane the region $D(r)$ lies 
within the rectangle $[-\sqrt{1-r^2}, \sqrt{1-r^2}]\times [-1, 0]$ 
and it is determined by the two parabolas $P_1: u=-r^2-p_r^2$ 
and $P_2: u=-\gamma-p_r^2$. See Figure 1 below 
that illustrates the two possible situations. 
\\[-15ex] 
\begin{center}
  \hspace{-10em}  
  \begin{minipage}{0.48\textwidth}
    \includegraphics[width=1.8\textwidth]{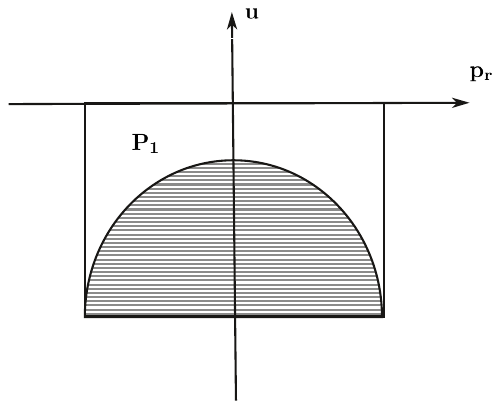}
  \end{minipage}
  \hspace{-5em} 
  \begin{minipage}{0.48\textwidth}
    \includegraphics[width=1.8\textwidth]{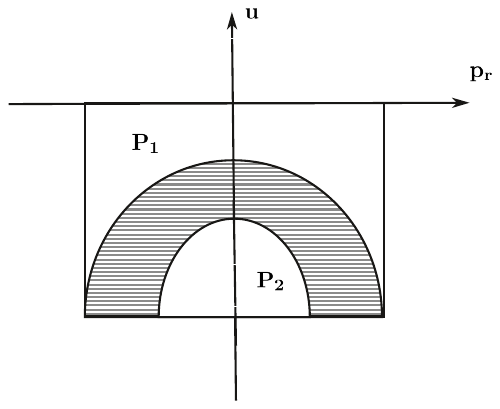}
  \end{minipage}
  \\[-65ex] 
  \captionof{figure}{$r^2\le\Gamma$ (left) and $r^2\ge\Gamma$ (right)}
\end{center}

%%%%%%%%%%%%%%%%%%%%%%%%%%%%%%%%%%%%%%%%%%%%%%%%%%%%%%%%%%%%%%%%%%%%%%%%%%%%%%%%%%

\setcounter{equation}{0}
\section{Construction of steady states}
\label{constr_sect} 

From the ansatz (\ref{fGam}) it is already possible 
to see that $f=f_\Gamma$ is a weak static solution, 
provided that we can prove later that 
\begin{equation}\label{rho_ziel} 
   \rho_f=\rho_{{\rm Kurth}}.
\end{equation} 
Let us verify this claim. 
Clearly $f(x, v)\ge 0$ for a.e.~$(x, v)\in\R^3\times\R^3$. 
If $f(x, v)\neq 0$, we must have 
\[ -1\le\ell^2-|x|^2-|v|^2\le 0,\quad\ell^2\le\Gamma. \] 
But the arguments in Section \ref{supp_sect} imply that these two relations 
enforce $|x|\le 1$ and $|v|\le 1$, cf.~(\ref{picc}). Therefore if we let 
\[ K=\{(x, v)\in\R^3\times\R^3: |x|\le 1, |v|\le 1\}, \] 
then $f(x, v)=0$ for a.e.~$(x, v)\in (\R^3\times\R^3)\setminus K$. 
According to (\ref{rho_ziel}), we also have $\rho_f\in L^\infty(\R^3)$. 

Next we are going to check that $f\in L^1(\R^3\times\R^3)$. 
To establish this assertion, we use Lemma \ref{intlem}(a) and Lemma \ref{Gamsuppc} 
on the support to write 
\begin{eqnarray}\label{freg1}   
   \lefteqn{\int_{\R^3}\int_{\R^3} f(x, v)\,dx\,dv}
   \nonumber \\ 
   & = & \int_{\R^3}\int_{\R^3}
   Q(\ell^2(x, v)-|x|^2-|v|^2)\,{\bf 1}_{\{\ell^2(x, v)\le\Gamma\}}\,dx\,dv
   \nonumber  
   \\ & = & \frac{1}{2}\int_0^\infty dr\int_{\R} dp_r\int_0^\infty d\ell^2\,
   Q(u)\,{\bf 1}_{\{\ell^2\le\Gamma\}} I_1(r, p_r, \ell)
   \nonumber     
   \\ & = & \frac{1}{2}\int_0^1 dr\,\frac{r^2}{1-r^2}\int\int_{D(r)}\,dp_r\,du 
   \,Q(u)\,I_1(r, p_r, \ell(r, p_r, u))
   \nonumber     
   \\ & = & \frac{1}{2}\int_0^{\sqrt{\Gamma}} dr\,\frac{r^2}{1-r^2}
   \int_{|p_r|\le\sqrt{1-r^2}} dp_r\int_{-1}^{-r^2-p_r^2} du 
   \,Q(u)\,I_1(r, p_r, \ell(r, p_r, u))
   \nonumber     
   \\ & & +\,\frac{1}{2}\int_{\sqrt{\Gamma}}^1 dr\,\frac{r^2}{1-r^2}
   \int_{|p_r|\le\sqrt{1-\gamma}} dp_r  
   \int_{-\gamma-p_r^2}^{-r^2-p_r^2} du
   \,Q(u)\,I_1(r, p_r, \ell(r, p_r, u))
   \nonumber     
   \\ & & +\,\frac{1}{2}\int_{\sqrt{\Gamma}}^1 dr\,\frac{r^2}{1-r^2}
   \int_{\sqrt{1-\gamma}\le |p_r|\le\sqrt{1-r^2}} dp_r  
   \int_{-1}^{-r^2-p_r^2} du\,Q(u)\,I_1(r, p_r, \ell(r, p_r, u)), 
   \qquad
\end{eqnarray} 
where for $\phi\in L^\infty(\R^3\times\R^3)$ in general 
\[ I_\phi(r, p_r, \ell)=\int_{{\rm SO}(3)}
   \phi\left(R\left(\begin{array}{c} 0 \\ 0 \\ r\end{array}\right),
   R\left(\begin{array}{c} \ell/r \\ 0 \\ p_r\end{array}\right)\right)\,dR, \] 
and we take $\phi=1$. See Lemma \ref{intlem}(a), 
and note that $\ell=\ell(r, p_r, u)\ge 0$ is determined from (\ref{ellu}).
Then in particular 
\begin{equation}\label{Iphi_bd} 
   |I_\phi(r, p_r, \ell)|\le\|\phi\|_{L^\infty(\R^3\times\R^3)}
   \int_{{\rm SO}(3)}\,dR=8\pi^2\|\phi\|_{L^\infty(\R^3\times\R^3)},
\end{equation}  
cf.~once again Lemma \ref{intlem}(a). Hence if we go back to (\ref{freg1}) 
and remember that $Q\ge 0$ is supposed to be integrable on $[-1, 0]$, 
we obtain the bound      
\begin{eqnarray*}   
   \lefteqn{\int_{\R^3}\int_{\R^3} f(x, v)\,dx\,dv}
   \\ & \le & 4\pi^2\int_0^{\sqrt{\Gamma}} dr\,\frac{r^2}{1-r^2}
   \int_{|p_r|\le\sqrt{1-r^2}} dp_r\int_{-1}^0 du\,Q(u)    
   \\ & & +\,4\pi^2\int_{\sqrt{\Gamma}}^1 dr\,\frac{r^2}{1-r^2}
   \int_{|p_r|\le\sqrt{1-\gamma}} dp_r  
   \int_{-1}^0 du\,Q(u)
   \\ & & +\,4\pi^2\int_{\sqrt{\Gamma}}^1 dr\,\frac{r^2}{1-r^2}
   \int_{\sqrt{1-\gamma}\le |p_r|\le\sqrt{1-r^2}} dp_r  
   \int_{-1}^0 du\,Q(u) 
   \\ & \le & 24\pi^2\,\|Q\|_{L^1([-1, 0])}\int_0^1\frac{r^2}{\sqrt{1-r^2}}\,dr. 
\end{eqnarray*} 
Since the right-hand side is finite, we conclude that 
$f\in L^1(\R^3\times\R^3)$. 

Lastly, to check that $f$ satisfies (\ref{distr_eq}), we are going to use 
the following general result on functions $f$ being defined 
in terms of conserved quantities. 

\begin{lemma}\label{inv_lem} Consider the autonomous system 
\begin{equation}\label{char_glen} 
   \dot{X}=V,\quad\dot{V}=-X,
\end{equation}  
in $\R^3\times\R^3$ and the associated flow 
\begin{equation}\label{flow_def} 
   \Lambda_t(x, v)=(X(t, x, v), V(t, x, v))
   =(x\cos t+v\sin t, -x\sin t+v\cos t).
\end{equation}  
Suppose that $f\in L^1(\R^3\times\R^3)$ is a function that is invariant 
under the flow, i.e., for every $t\in\R$ we have 
\begin{equation}\label{fphia} 
   f(\Lambda_t(x, v))=f(x, v)\quad\mbox{for a.e.}\quad (x, v)\in\R^3\times\R^3.
\end{equation}  
If the (essential) support of $f$ is contained 
in the set $\{(x, v)\in\R^3\times\R^3: |x|\le 1\}$ and if 
\begin{equation}\label{rho1} 
   \rho_f(x)=\rho_{{\rm Kurth}}(x)=\frac{3}{4\pi}
   \,{\bf 1}_{\overline{B_1(0)}}\,(x)\quad\mbox{for a.e.}\quad x\in\R^3,
\end{equation} 
then $f$ solves (\ref{distr_eq}). 
\end{lemma} 
{\bf Proof of Lemma \ref{inv_lem}:} First note that (\ref{rho1}) leads to 
\[ U_f(x)=-\int_{\R^3}\frac{\rho_f(y)}{|x-y|}\,dy
   =\left\{\begin{array}{c@{\quad:\quad}c} 
   \frac{1}{2} |x|^2-\frac{3}{2} & |x|\le 1
   \\[1ex] -\frac{1}{|x|} & |x|>1
   \end{array}\right. \] 
by direct calculation. In particular, we obtain $\nabla_x U_f(x)=x$ 
for $|x|<1$. Due to the support assumption on $f$ we hence need to verify that 
\[ A:=\int_{\R^3}\int_{\R^3} f(x, v)
   \,\Big(v\cdot\nabla_x\varphi(x, v)-x\cdot\nabla_v\varphi(x, v)\Big)\,dx\,dv=0 \] 
for $\varphi\in {\cal D}(\R^3\times\R^3)$. To begin with, 
since $\Lambda_{2\pi}=\Lambda_0$, we observe that by (\ref{char_glen}) 
for all $(x, v)\in\R^3\times\R^3$ 
\begin{eqnarray}\label{02pi}  
   & & \int_0^{2\pi}\Big[V(t, x, v)\cdot\nabla_x\varphi(\Lambda_t(x, v))
   -X(t, x, v)\cdot\nabla_v\varphi(\Lambda_t(x, v))\Big]\,dt 
   \nonumber   
   \\ & & =\,\int_0^{2\pi} \frac{d}{dt}\,\varphi(\Lambda_t(x, v))\,dt 
   \nonumber   
   \\ & & =\,\varphi(\Lambda_{2\pi}(x, v))-\varphi(\Lambda_0(x, v))
   =0.  
\end{eqnarray} 
Using that $\Lambda_t$ is measure-preserving, we can apply (\ref{fphia}) 
for $t\in\R$ to get 
\begin{eqnarray*} 
   A & = & \int_{\R^3}\int_{\R^3} f(\Lambda_t(x, v))
   \,\Big(V(t, x, v)\cdot\nabla_x\varphi(\Lambda_t(x, v))
   -X(t, x, v)\cdot\nabla_v\varphi(\Lambda_t(x, v))\Big)\,dx\,dv
   \\ & = & \int_{\R^3}\int_{\R^3} f(x, v)
   \,\Big(V(t, x, v)\cdot\nabla_x\varphi(\Lambda_t(x, v))
   -X(t, x, v)\cdot\nabla_v\varphi(\Lambda_t(x, v))\Big)\,dx\,dv. 
\end{eqnarray*} 
Integrating this relation on $[0, 2\pi]$ in $t$ and observing 
(\ref{02pi}), it follows that $2\pi A=0$, and hence $A=0$. 
{\hfill$\Box$}\bigskip 
 
Since both $|x|^2+|v|^2$ and $|x\wedge v|^2$ 
are invariant under the flow (\ref{flow_def}), 
for $f$ as in (\ref{fGam}) it follows that 
$f$ is a weak static solution. 
Hence in order to complete the proof of Theorem \ref{ex_fGam}, 
it remains to satisfy the condition (\ref{rho_ziel}).  
\smallskip 

This question will be approached in the next section, 
but before we need to record a representation formula 
that will be useful later. Exactly as in (\ref{freg1}) 
it is shown that given $\phi\in L^\infty(\R^3\times\R^3)$ one has 
\begin{eqnarray}\label{freg}   
   \lefteqn{\int_{\R^3}\int_{\R^3} f(x, v)\,\phi(x, v)\,dx\,dv}
   \nonumber 
   \\ & = & \frac{1}{2}\int_0^{\sqrt{\Gamma}} dr\,\frac{r^2}{1-r^2}
   \int_{|p_r|\le\sqrt{1-r^2}} dp_r\int_{-1}^{-r^2-p_r^2} du 
   \,Q(u)\,I_\phi(r, p_r, \ell(r, p_r, u))
   \nonumber     
   \\ & & +\,\frac{1}{2}\int_{\sqrt{\Gamma}}^1 dr\,\frac{r^2}{1-r^2}
   \int_{|p_r|\le\sqrt{1-\gamma}} dp_r  
   \int_{-\gamma-p_r^2}^{-r^2-p_r^2} du
   \,Q(u)\,I_\phi(r, p_r, \ell(r, p_r, u))
   \nonumber     
   \\ & & +\,\frac{1}{2}\int_{\sqrt{\Gamma}}^1 dr\,\frac{r^2}{1-r^2}
   \int_{\sqrt{1-\gamma}\le |p_r|\le\sqrt{1-r^2}} dp_r  
   \int_{-1}^{-r^2-p_r^2} du\,Q(u)\,I_\phi(r, p_r, \ell(r, p_r, u)).  
   \qquad
\end{eqnarray} 

%%%%%%%%%%%%%%%%%%%%%%%%%%%%%%%%%%%%%%%%%%%%%%%%%%%%%%%%%%%%%%%%%%%%%%%%%%%%%%%%%%

\setcounter{equation}{0}
\section{Reformulation of the condition on uniform distribution}
\label{reform_sect} 

Let us specify the ansatz (\ref{fGam}) further. 
Suppose that $\varphi: [0, 1]\to\R$ is absolutely continuous, 
monotone non-decreasing and such that $\varphi(0)=0$. 
In other words, $\varphi\in W^{1, 1}(]0, 1[)$, the Sobolev space, 
and $\varphi(0)=0$ as well as $\varphi'(s)\ge 0$ a.e. 
Define $\Phi(u)=\varphi(u+1)$ for $u\in [-1, 0]$ and $\Phi(u)=0$ 
for $u\in\R\setminus [-1, 0]$. Lastly, let 
\begin{equation}\label{konsu} 
   Q(u)=\left\{\begin{array}{c@{\quad:\quad}c} \frac{3}{4\pi}\,\Phi'(u)
   & {\rm a.e.}\,\,u\in ]0, 1[ 
   \\[1ex] 0 & {\rm else}\end{array}\right. .
\end{equation}  
Note that all those functions will depend on $\Gamma$, 
and the factor $\frac{3}{4\pi}$ is introduced to match 
the Kurth solution. Then $Q=Q(u): \R\to [0, \infty[$ 
is an integrable function with support in $[-1, 0]$. 

For the density $\rho_f=\rho_{f_\Gamma}$ 
induced by (\ref{fGam}) and given by $\rho_f(x)=\int_{\R^3} f(x, v)\,dv$, 
we first observe that $\rho_f(x)=0$ for $|x|>1$ by the support properties of $f$. 
If $r=|x|\le 1$, then Lemma \ref{intlem}(b) implies
\begin{eqnarray*} 
   \rho_f(x) & = & \frac{3}{4\pi}\int_{\R^3}\Phi'(\ell^2(x, v)-|x|^2-|v|^2)
   \,{\bf 1}_{\{\ell^2(x, v)\le\Gamma\}}\,dv 
   \\ & = & \frac{3}{4r^2}\int_{\R} dp_r\int_0^\infty d\ell^2
   \,\Phi'\Big(-r^2-p_r^2-\frac{\ell^2}{r^2}(1-r^2)\Big)\,{\bf 1}_{\{\ell^2\le\Gamma\}}
   \\ & = & \frac{3}{4(1-r^2)}\int_{\R} dp_r\int_{-1}^0 du
   \,\Phi'(u)\,{\bf 1}_{\{\ell^2(r, p_r, u)\le\Gamma\}}\,. 
\end{eqnarray*} 
In Lemma \ref{Gamsuppc} we have determined the support properties of the integrand. 
\smallskip 

\noindent 
(i) If $0\le r\le\sqrt{\Gamma}$, then Lemma \ref{Gamsuppc}(b) together with $\Phi(-1)=0$ yields 
\begin{eqnarray*} 
   \rho_f(x) & = & \frac{3}{4(1-r^2)}
   \int_{|p_r|\le\sqrt{1-r^2}} dp_r\int_{-1}^{-r^2-p_r^2} du\,\Phi'(u)
   \\ & = & \frac{3}{2(1-r^2)}\int_0^{\sqrt{1-r^2}} \Phi(-r^2-p_r^2)\,dp_r
   \\ & = & \frac{3}{4\pi}\frac{\pi}{1-r^2}\int_0^{1-r^2} \Phi(-(1-y))\,\frac{dy}{\sqrt{(1-r^2)-y}}. 
\end{eqnarray*}

\noindent 
(ii) If $\sqrt{\Gamma}\le r<1$, then by Lemma \ref{Gamsuppc}(c) we obtain 
\begin{eqnarray*} 
   \rho_f(x) 
   & = & \frac{3}{4(1-r^2)}\int_{|p_r|\le\sqrt{1-\gamma}} dp_r
   \int_{-\gamma-p_r^2}^{-r^2-p_r^2} du\,\Phi'(u)
   \\ & & +\,\frac{3}{4(1-r^2)}\int_{\sqrt{1-\gamma}\le |p_r|\le\sqrt{1-r^2}} dp_r
   \int_{-1}^{-r^2-p_r^2} du\,\Phi'(u)
   \\ & = & \frac{3}{2(1-r^2)}\int_0^{\sqrt{1-\gamma}}
   \Big[\Phi(-r^2-p_r^2)-\Phi(-\gamma-p_r^2)\Big]\,dp_r
   \\ & & +\,\frac{3}{2(1-r^2)}\int_{\sqrt{1-\gamma}}^{\sqrt{1-r^2}}
   \,\Phi(-r^2-p_r^2)\,dp_r
   \\ & = & \frac{3}{2(1-r^2)}\int_0^{\sqrt{1-r^2}}\Phi(-r^2-p_r^2)\,dp_r
   -\frac{3}{2(1-r^2)}\int_0^{\sqrt{1-\gamma}}
   \Phi\Big(-\Big(\frac{1-r^2}{r^2}\Big)\Gamma-r^2-p_r^2\Big)\,dp_r
   \\ & = & \frac{3}{4\pi}\,\frac{\pi}{1-r^2}\bigg[ 
   \int_0^{1-r^2}\Phi(-(1-y))\,\frac{dy}{\sqrt{(1-r^2)-y}}
   \\ & & \hspace{4em} -\,\int^{1-r^2}_{\Gamma(\frac{1-r^2}{r^2})}
   \Phi\Big(-\Gamma\Big(\frac{1-r^2}{r^2}\Big)-(1-y)\Big)\,\frac{dy}{\sqrt{(1-r^2)-y}}\bigg].  
\end{eqnarray*}
To summarize, under the assumptions made above on $\varphi(s)$ 
or the shifted $\Phi(u)=\varphi(u+1)$, 
(\ref{rho_ziel}) will hold if and only if 
\begin{eqnarray}\label{rhoeq2} 
   1=\left\{\begin{array}{l@{\quad}l} 
   \frac{\pi}{1-r^2}\int_0^{1-r^2}\Phi(-(1-y))
   \,\frac{dy}{\sqrt{(1-r^2)-y}} & :\quad 0\le r\le\sqrt{\Gamma}
   \\[1ex] \frac{\pi}{1-r^2}\bigg[ 
   \int_0^{1-r^2}\Phi(-(1-y))\,\frac{dy}{\sqrt{(1-r^2)-y}}
   \\ \hspace{3em} -\,\int^{1-r^2}_{\Gamma(\frac{1-r^2}{r^2})}
   \Phi(-\Gamma(\frac{1-r^2}{r^2})-(1-y))\,\frac{dy}{\sqrt{(1-r^2)-y}}\bigg] 
   & :\quad\sqrt{\Gamma}\le r<1\end{array}\right. , 
\end{eqnarray} 
is verified. In the following Sections \ref{prel_sect}-\ref{PhiGam_sect}, 
(\ref{rhoeq2}) will be rewritten and solved 
as a suitable Abel-type integral equation. 

%%%%%%%%%%%%%%%%%%%%%%%%%%%%%%%%%%%%%%%%%%%%%%%%%%%%%%%%%%%%%%%%%%%%%%%%%%%%%%%%%%%%%

\setcounter{equation}{0}
\section{Some integral operators}
\label{prel_sect} 

In this section we will mostly be dealing with functions in $C([0, 1])$. Let 
\[ (I^{1/2}g)(s)=\int_0^s\frac{g(\sigma)}{\sqrt{s-\sigma}}\,d\sigma
   =\int_0^s\frac{g(s-\sigma)}{\sqrt{\sigma}}\,d\sigma,
   \quad s\in [0, 1], \] 
and, for $0<\Gamma<1$,  
\[ \chi_\Gamma(s)=\left\{\begin{array}{c@{\quad:\quad}c} 
   \Gamma\frac{s}{1-s} & s\in [0, 1-\Gamma]
   \\[1ex] s & s\in [1-\Gamma, 1]\end{array}\right. . \]
Furthermore, we define 
\[ (K_\Gamma g)(s)=\int_{\chi_\Gamma(s)}^s 
   g(\sigma-\chi_\Gamma(s))\,\frac{d\sigma}{\sqrt{s-\sigma}}
   =(I^{1/2}g)(s-\chi_\Gamma(s)),
   \quad s\in [0, 1]. \] 
Thus if we let 
\[ (S_\Gamma\psi)(s)=\psi(s-\chi_\Gamma(s)), \] 
we have
\begin{equation}\label{ento} 
   K_\Gamma g=S_\Gamma I^{1/2}g.
\end{equation}  
As it will turn out later, in order to obtain (\ref{rhoeq2}) 
we need to solve the equation 
\[ (I^{1/2}\varphi)(s)-(K_\Gamma\varphi)(s)=\frac{s}{\pi},
   \quad s\in [0, 1], \] 
for $\varphi$, which, due to (\ref{ento}), can be written as 
\[ (({\rm id}-S_\Gamma)I^{1/2}\varphi)(s)=\frac{s}{\pi}; \]    
here ${\rm id}$ denotes the identity operator. 
Thus we first have to study the equation $(({\rm id}-S_\Gamma)\psi)(s)=s$ 
for $s\in [0, 1]$. 

%%%%%%%%%%%%%%%%%%%%%%%%%%%%%%%%%%%%%%%%%%%%%%%%%%%%%%%%%%%%%%%%%%%%%%%%%%%%%%%%%%%

\setcounter{equation}{0}
\section{Solving the functional equation 
$(({\rm id}-S_\Gamma)\psi)(s)=s$}

We recall that $(S_\Gamma\psi)(s)=\psi(h_\Gamma(s))$ 
for $h_\Gamma(s)=s-\chi_\Gamma(s)$ and $s\in [0, 1]$. 
To simplify notation we will write $h$ instead of $h_\Gamma$ 
throughout this section. 
By $h^k=h\circ h^{k-1}$ we denote the iterates of 
\begin{equation}\label{hdef} 
   h(s)=\left\{\begin{array}{c@{\quad:\quad}c} 
   s(1-\Gamma\frac{1}{1-s}) & s\in [0, 1-\Gamma]
   \\[1ex] 0 & s\in [1-\Gamma, 1]\end{array}\right. ,
\end{equation} 
where $h^0(s)=s$.       
   
\begin{lemma}\label{le1} 
Let $\Gamma\in [\frac{63}{64}, 1[$ and define $R=2(1-\Gamma)^{1/3}<1$ 
as well as $h_c(z)=z(1-\Gamma\,\frac{1}{1-z})$ for $z\in\C$ 
such that $|z|\le R$. Then $|h_c(z)|\le 6(1-\Gamma)^{1/3}|z|$ for $|z|\le R$. 
\end{lemma} 
{\bf Proof\,:} In the domain $|z|\le R$ we have 
\[ \Big|1-\Gamma\,\frac{1}{1-z}\Big|
   =\Big|\frac{1-z-\Gamma}{1-z}\Big|
   \le\Big|\frac{-R-(1-\Gamma)}{R-1}\Big|
   \le\frac{3(1-\Gamma)^{1/3}}{1-R}\le 6(1-\Gamma)^{1/3}, \] 
noting that $R\le 1/2$.  
{\hfill$\Box$}\bigskip 

\begin{lemma}\label{le2} 
Let $\Gamma\in [\frac{1727}{1728}, 1[$ 
and define $R=2(1-\Gamma)^{1/3}$ 
as well as $h_c(z)=z(1-\Gamma\,\frac{1}{1-z})$ for $z\in\C$ 
such that $|z|\le R$. 
\begin{itemize} 
\item[(a)] The series 
\begin{equation}\label{psitil} 
   \tilde{\psi}(z)=\sum_{k=0}^\infty h_c^k(z),
   \quad |z|<R,
\end{equation}  
is well-defined and defines a holomorphic function such that 
\begin{equation}\label{functequ} 
   \tilde{\psi}(z)-\tilde{\psi}(h_c(z))=z,\quad |z|<R.
\end{equation} 
Moreover, $\tilde{\psi}(0)=0$, $\tilde{\psi}(1-\Gamma)=1-\Gamma$, 
$\tilde{\psi}'(0)=\frac{1}{\Gamma}$ and 
$\tilde{\psi}'(1-\Gamma)=1-\frac{1}{\Gamma^2}(1-\Gamma)$. 
\item[(b)] We have the following estimates: 
\begin{eqnarray}
   & & \sup_{|z|<R}|\tilde{\psi}(z)|\le 4(1-\Gamma)^{1/3}; 
   \nonumber    
   \\ & & \sup_{|z|\le 1-\Gamma} |\tilde{\psi}'(z)|\le 16;
   \label{1abl}    
   \\ & & \sup_{|z|\le 1-\Gamma} |\tilde{\psi}''(z)|
   \le\frac{64}{(1-\Gamma)^{1/3}}.   
   \label{2abl} 
\end{eqnarray}  
\end{itemize} 
\end{lemma}
{\bf Proof\,:} (a) Using Lemma \ref{le1} and the restriction on $\Gamma$, 
we get $|h_c(z)|\le 6(1-\Gamma)^{1/3} |z|\le (1/2)|z|$ for $|z|\le R$.
Therefore $|h_c^k(z)|\le 2^{-k}|z|\le 2^{-k+1}(1-\Gamma)^{1/3}$ 
shows that (\ref{psitil}) is a series of uniformly convergent 
holomorphic functions. In particular, we obtain 
\[ \tilde{\psi}(z)-\tilde{\psi}(h_c(z))
   =\sum_{k=0}^\infty h_c^k(z)
   -\sum_{k=0}^\infty h_c^{k+1}(z)=h_c^0(z)=z,
   \quad |z|<R. \] 
Since $h_c(0)=0$, it follows that $\tilde{\psi}(0)=0$. 
Differentiating (\ref{functequ}) and observing $h'_c(z)=1-\frac{\Gamma}{(1-z)^2}$, 
the relation $\Gamma\tilde{\psi}'(0)=[1-h'_c(0)]\tilde{\psi}'(0)=1$ is found, 
which gives $\tilde{\psi}'(0)=\frac{1}{\Gamma}$. 
In addition, $h_c(1-\Gamma)=0$ by definition, and hence $h_c(0)=0$ 
yields $h_c^k(1-\Gamma)=0$ for $k\in\N$. As a consequence, 
we must have $\tilde{\psi}(1-\Gamma)=1-\Gamma$. 
Similarly, $1=\tilde{\psi}'(1-\Gamma)-\tilde{\psi}'(h_c(1-\Gamma))
h'_c(1-\Gamma)=\tilde{\psi}'(1-\Gamma)-\tilde{\psi}'(0)(1-\frac{1}{\Gamma})
=\tilde{\psi}'(1-\Gamma)-\frac{1}{\Gamma}(1-\frac{1}{\Gamma})$ 
leads to $\tilde{\psi}'(1-\Gamma)=1-\frac{1}{\Gamma^2}(1-\Gamma)$. 
(b) To bound $\tilde{\psi}$, we can use (a) and estimate
\[ |\tilde{\psi}(z)|\le\sum_{k=0}^\infty |h_c^k(z)|
   \le\sum_{k=0}^\infty 2^{-k+1}(1-\Gamma)^{1/3}=4(1-\Gamma)^{1/3}. \] 
For the derivatives, we rely on the Cauchy estimates \cite[Satz 5.1]{FL} 
for general functions $f$ that are holomorphic 
in an open neighborhood of $B_r(0)$, $\delta\in ]0, r]$, 
$j\in\N_0$, $|z|\le r-\delta$: then $|f^{(j)}(z)|\le\frac{r}{\delta}\,\frac{j!}{\delta^j}
\,\max_{|\zeta|=r} |f(\zeta)|$. 
We apply this bound with $r=(1-\Gamma)^{1/3}$ 
and $\delta=(1-\Gamma)^{1/3}-(1-\Gamma)
=(1-\Gamma)^{1/3}[1-(1-\Gamma)^{2/3}]
\ge (1/2)(1-\Gamma)^{1/3}$. Therefore 
\[ |\tilde{\psi}'(z)|\le\frac{4}{(1-\Gamma)^{1/3}}
   \,\max_{|\zeta|=(1-\Gamma)^{1/3}} |\tilde{\psi}(\zeta)|
   \le 16 \] 
and 
\[ |\tilde{\psi}''(z)|\le 
   \,\frac{16}{(1-\Gamma)^{2/3}}
   \,\max_{|\zeta|=(1-\Gamma)^{1/3}} |\tilde{\psi}(\zeta)|
   \le\frac{64}{(1-\Gamma)^{1/3}} \] 
are obtained, both for $|z|\le 1-\Gamma$; note that $r-\delta=1-\Gamma$. 
{\hfill$\Box$}\bigskip 

\begin{lemma}\label{le3} In the setup of Lemma \ref{le2}, define 
\begin{equation}\label{psidef} 
   \psi(s)=\left\{\begin{array}{c@{\quad:\quad}c} 
   \tilde{\psi}(s) & s\in [0, 1-\Gamma]
   \\[1ex] s & s\in [1-\Gamma, 1]\end{array}\right. .
\end{equation} 
Then $\psi: [0, 1]\to\R$ is piecewise smooth, 
in the sense that $\psi$ is continuous and its restrictions 
to $[0, 1-\Gamma]$ and $[1-\Gamma, 1]$ are $C^\infty$; 
in particular, $\psi$ is Lipschitz continuous. Moreover, 
\begin{equation}\label{functequx} 
   \psi(s)-\psi(h(s))=s,\quad s\in [0, 1],
\end{equation} 
as well as $\psi(0)=0$ and $\psi'(0)=\frac{1}{\Gamma}$. 
\end{lemma}
{\bf Proof\,:} First observe that $1-\Gamma<R=2(1-\Gamma)^{1/3}$, 
which means that $\tilde{\psi}(s)$ is defined for $s\in [0, 1-\Gamma]$. 
Next note that $\psi$ is continuous at $s=1-\Gamma$, 
as $\tilde{\psi}(1-\Gamma)=1-\Gamma$ by Lemma \ref{le2}(a). 
Also $\psi$ is $C^\infty$ in $[0, 1-\Gamma]$ 
with 
\[ \psi'((1-\Gamma)^-)=1-\frac{1}{\Gamma^2}(1-\Gamma), \]  
and furthermore $\psi$ is $C^\infty$ in $[1-\Gamma, 1]$ 
such that 
\[ \psi'((1-\Gamma)^+)=1 \] 
does exist; higher-order derivatives can be obtained inductively, 
To establish (\ref{functequx}), first let $s\in [0, 1-\Gamma]$.  
Then (\ref{functequ}) shows that (\ref{functequx}) is verified. 
If $s\in [1-\Gamma, 1]$, then $h(s)=0$ by (\ref{hdef}), 
and thus $\psi(h(s))=\psi(0)=\tilde{\psi}(0)=0$ 
due to Lemma \ref{le2}. It follows that (\ref{functequx}) 
holds for all $s\in [0, 1]$. In addition, $\psi(0)=\tilde{\psi}(0)=0$ 
and $1=\psi'(0)-\psi'(h(0))h'(0)=\psi'(0)-\psi'(0)(1-\Gamma)=\Gamma\psi'(0)$ 
by Lemma \ref{le2}(a) and (\ref{functequx}). 
{\hfill$\Box$}\bigskip 

\begin{remark}\label{Gam1_1}  
{\rm The construction of $\psi$ extends to the case $\Gamma=1$ as follows. 
Firstly, we get $h(s)=0$ for $s\in [0, 1]$ in (\ref{hdef}). 
By (\ref{psidef}) we obtain $\psi(s)=s$ for $s\in [0, 1]$.
}  
\end{remark} 

%%%%%%%%%%%%%%%%%%%%%%%%%%%%%%%%%%%%%%%%%%%%%%%%%%%%%%%%%%%%%%%%%%%%%%%%%%%%%%%%%%%

\setcounter{equation}{0}
\section{Solving the integral equation 
$(({\rm id}-S_\Gamma)I^{1/2}\varphi)(s)=\frac{s}{\pi}$}

\begin{lemma}\label{le4} 
For $\Gamma\in [\frac{1727}{1728}, 1[$
let $\psi=\psi_\Gamma$ be given by (\ref{psidef}) in Lemma \ref{le3}.  
Define $\varphi=\varphi_\Gamma$ by means of 
\begin{equation}\label{varphiform}  
   \varphi(s)=\frac{1}{\pi^2}
   \frac{d}{ds}\int_0^s\frac{\psi(s-\sigma)}{\sqrt{\sigma}}\,d\sigma,
   \quad s\in ]0, 1].
\end{equation}  
Then $\varphi\in\bigcap_{\,0<\beta<\frac{1}{2}} C^{0, \beta}([0, 1])$, 
$\varphi(0)=0$, and 
\begin{equation}\label{zielgl} 
   (({\rm id}-S_\Gamma)I^{1/2}\varphi)(s)=\frac{s}{\pi},
   \quad s\in [0, 1].
\end{equation}  
\end{lemma} 
{\bf Proof\,:} The equation (\ref{zielgl}) reads as  
\[ (I^{1/2}\varphi)(s)-(I^{1/2}\varphi)(h(s))
   =\frac{s}{\pi},\quad s\in [0, 1], \] 
so according to (\ref{functequx}) in Lemma \ref{le3}, 
it suffices to establish that 
\begin{equation}\label{I12} 
   \pi\,(I^{1/2}\varphi)(s)=\psi(s),\quad s\in [0, 1].
\end{equation}  
For this, let 
\[ (J^{1/2}g)(s)=\frac{1}{\pi}\int_0^s\frac{g(\sigma)}{\sqrt{s-\sigma}}\,d\sigma
   =\frac{1}{\pi}\,(I^{1/2}g)(s),
   \quad s\in [0, 1], \] 
denote the standard Abel integral operator, 
see \cite{HL} and \cite[equ.~(1.1.1)]{GV}. 
Then (\ref{I12}) reads as 
\begin{equation}\label{J12} 
   (J^{1/2}\varphi)(s)=\frac{1}{\pi^2}\,\psi(s),\quad s\in [0, 1].
\end{equation}  
From Lemma \ref{le3} we know that $\psi: [0, 1]\to\R$ is Lipschitz 
continuous and such that $\psi(0)=0$. In particular, 
$\psi\in C^{0, \frac{1}{2}+\beta}([0, 1])$, 
the H\"older space, for every $0<\beta<\frac{1}{2}$. 
According to \cite[Thm.~5.1.1]{GV}, (\ref{J12}) has a unique solution 
$\varphi\in C^{0, \beta}([0, 1])$, and $\varphi(0)=0$. 
In addition, the solution is given by the formula 
\[ \varphi(s)=\frac{1}{\pi^2}
   \frac{d}{ds}\int_0^s\frac{\psi(\sigma)\,d\sigma}{\sqrt{s-\sigma}},
   \quad s\in [0, 1]. \] 
This completes the proof of the lemma. 
{\hfill$\Box$}\bigskip 

\begin{lemma}\label{le5} In the setup of Lemma \ref{le4}, 
suppose that in addition $\Gamma\in [1-2^{-12}, 1[$ holds. Then 
\begin{itemize} 
\item[(a)] $\varphi\in C^\infty(]0, 1-\Gamma[)$ and 
\[ \frac{1}{2\Gamma\pi^2\sqrt{s}}\le\varphi'(s)
   \le\frac{1}{\pi^2}\,\frac{1}{\Gamma\sqrt{s}}
   +\frac{1}{\pi^2}\frac{128}{(1-\Gamma)^{1/3}}\,\sqrt{s},
   \quad s\in ]0, 1-\Gamma[. \]  
\item[(b)] $\varphi\in C^\infty(]1-\Gamma, 1])$ and 
\[ \frac{1}{\Gamma^2\pi^2}\frac{1-\Gamma}{\sqrt{s-(1-\Gamma)}}
   \le\varphi'(s)
   \le\frac{2+128\,(1-\Gamma)^{2/3}}{\pi^2\sqrt{s}}
   +\frac{1}{\Gamma^2\pi^2}\frac{1-\Gamma}{\sqrt{s-(1-\Gamma)}}, 
   \quad s\in ]1-\Gamma, 1]. \]    
\end{itemize}  
\end{lemma} 
{\bf Proof\,:} (a) We consider $s\in [0, 1-\Gamma[$. 
Since $\psi$ is $C^\infty$ on this interval by Lemma \ref{le3} 
and $\psi(0)=0$, we obtain   
\[ \varphi(s)=\frac{1}{\pi^2}
   \int_0^s\frac{\psi'(s-\xi)}{\sqrt{\xi}}\,d\xi \] 
as well as 
\begin{equation}\label{mudy} 
   \varphi'(s)=\frac{1}{\pi^2}\,\frac{\psi'(0)}{\sqrt{s}}
   +\frac{1}{\pi^2}\int_0^s\frac{\psi''(s-\sigma)}{\sqrt{\sigma}}\,d\sigma
   =\frac{1}{\pi^2}\,\frac{1}{\Gamma\sqrt{s}}
   +\frac{1}{\pi^2}\int_0^s\frac{\psi''(s-\sigma)}{\sqrt{\sigma}}\,d\sigma,
\end{equation}  
cf.~Lemma \ref{le3}. For the second term, (\ref{2abl}) yields 
\begin{equation}\label{2ndt} 
    \Big|\int_0^s\frac{\psi''(s-\sigma)}{\sqrt{\sigma}}\,d\sigma\Big|
   \le\frac{64}{(1-\Gamma)^{1/3}}\int_0^s\frac{1}{\sqrt{\sigma}}\,d\sigma
   =\frac{128}{(1-\Gamma)^{1/3}}\,\sqrt{s}.
\end{equation}  
Therefore 
\begin{eqnarray*} 
   \varphi'(s) 
   & \ge & \frac{1}{\pi^2}\,\frac{1}{\Gamma\sqrt{s}}
   -\frac{1}{\pi^2}\,\frac{128}{(1-\Gamma)^{1/3}}\,\sqrt{s}
   \\ & = & \frac{1}{\pi^2\sqrt{s}}\Big(\frac{1}{\Gamma}
   -\frac{128}{(1-\Gamma)^{1/3}}\,s\Big)
   \\ & \ge & \frac{1}{\pi^2\sqrt{s}}\Big(\frac{1}{\Gamma}
   -128(1-\Gamma)^{2/3}\Big)
   \\ & \ge & \frac{1}{2\Gamma\pi^2\sqrt{s}}, 
\end{eqnarray*} 
observing the additional restriction on $\Gamma$ and $\frac{1}{2\Gamma}\ge\frac{1}{2}$. 
The upper bound follows similarly. 
The higher order derivatives of $\varphi$ can be obtained 
by taking further derivatives of (\ref{mudy}). 
(b) For $s\in ]1-\Gamma, 1]$ we can rewrite (\ref{varphiform}) as 
\begin{eqnarray*} 
   \varphi(s) & = & \frac{1}{\pi^2}
   \frac{d}{ds}\int_0^s\frac{\psi(\sigma)}{\sqrt{s-\sigma}}\,d\sigma
   \\ & = & \frac{1}{\pi^2}
   \frac{d}{ds}\Big(\int_0^{1-\Gamma}\frac{\psi(\sigma)}{\sqrt{s-\sigma}}\,d\sigma
   +\int_{1-\Gamma}^s\frac{\sigma}{\sqrt{s-\sigma}}\,d\sigma\Big) 
   \\ & = & \frac{1}{\pi^2}\frac{d}{ds}
   \int_{s-(1-\Gamma)}^s\frac{\psi(s-\sigma)}{\sqrt{\sigma}}\,d\sigma
   +\frac{1}{\pi^2}\frac{d}{ds}
   \bigg[s^{3/2}\bigg(-\frac{2}{3}\Big(1-\frac{1-\Gamma}{s}\Big)^{3/2}
   +2\sqrt{1-\frac{1-\Gamma}{s}}\bigg)\bigg]
   \\ & = & \frac{1}{\pi^2}\frac{\psi(0)}{\sqrt{s}}
   -\frac{1}{\pi^2}\frac{\psi(1-\Gamma)}{\sqrt{s-(1-\Gamma)}}
   +\frac{1}{\pi^2}
   \int_{s-(1-\Gamma)}^s\frac{\psi'(x-\sigma)}{\sqrt{\sigma}}\,d\sigma
   \\ & & +\,\frac{1}{\pi^2}
   \bigg[\frac{s}{\sqrt{s-(1-\Gamma)}}
   +\sqrt{s-(1-\Gamma)}\bigg]
   \\ & = & -\frac{1}{\pi^2}\frac{(1-\Gamma)}{\sqrt{s-(1-\Gamma)}}
   +\frac{1}{\pi^2}\int_{s-(1-\Gamma)}^s\frac{\psi'(s-\sigma)}{\sqrt{\sigma}}\,d\sigma
   \\ & & +\,\frac{1}{\pi^2}
   \bigg[\frac{s}{\sqrt{s-(1-\Gamma)}}
   +\sqrt{s-(1-\Gamma)}\bigg]
   \\ & = & \frac{1}{\pi^2}\int_{s-(1-\Gamma)}^s\frac{\psi'(s-\sigma)}{\sqrt{\sigma}}\,d\sigma
   +\frac{2}{\pi^2}\sqrt{s-(1-\Gamma)}. 
\end{eqnarray*} 
As a consequence, 
\begin{eqnarray}\label{2ablb}  
   \varphi'(s) 
   & = & \frac{1}{\pi^2}\frac{\psi'(0)}{\sqrt{s}}
   -\frac{1}{\pi^2}\frac{\psi'((1-\Gamma)^-)}{\sqrt{s-(1-\Gamma)}}
   +\frac{1}{\pi^2}\int_{s-(1-\Gamma)}^s\frac{\psi''(s-\sigma)}{\sqrt{\sigma}}\,d\sigma
   +\frac{1}{\pi^2}\frac{1}{\sqrt{s-(1-\Gamma)}} 
   \nonumber   
   \\ & = & \frac{1}{\pi^2}\frac{1}{\Gamma\sqrt{s}}
   +\frac{1}{\pi^2}\int_{s-(1-\Gamma)}^s\frac{\psi''(s-\sigma)}{\sqrt{\sigma}}\,d\sigma
   +\frac{1}{\pi^2}\frac{1-\Gamma}{\Gamma^2\sqrt{s-(1-\Gamma)}},   
\end{eqnarray} 
cf.~Lemma \ref{le3} and its proof. 
Now we can apply (\ref{2abl}) once again to deduce that  
\begin{eqnarray}\label{flopl}   
   \Big|\int_{s-(1-\Gamma)}^s\frac{\psi''(s-\sigma)}{\sqrt{\sigma}}\,d\sigma\Big|
   & \le & \frac{64}{(1-\Gamma)^{1/3}}\int_{s-(1-\Gamma)}^s
   \frac{1}{\sqrt{\sigma}}\,d\sigma
   \nonumber
   \\ & = & 
   \frac{128}{(1-\Gamma)^{1/3}}\,\Big(\sqrt{s}-\sqrt{s-(1-\Gamma)}\Big)
   \nonumber   
   \\ & = & 
   128\,(1-\Gamma)^{2/3}\,\frac{1}{\sqrt{s}+\sqrt{s-(1-\Gamma)}}
   \nonumber   
   \\ & \le & 128\,(1-\Gamma)^{2/3}\,\frac{1}{\sqrt{s}}.  
\end{eqnarray} 
From the condition on $\Gamma$, it follows that 
\begin{eqnarray*} 
   \varphi'(s) 
   & \ge & \frac{1}{\pi^2}\frac{1}{\Gamma\sqrt{s}}
   -\frac{128}{\pi^2}\,(1-\Gamma)^{2/3}\,\frac{1}{\sqrt{s}}
   +\frac{1}{\pi^2}\frac{1-\Gamma}{\Gamma^2\sqrt{s-(1-\Gamma)}}
   \\ & \ge & \frac{1}{\pi^2}\frac{1-\Gamma}{\Gamma^2\sqrt{s-(1-\Gamma)}},    
\end{eqnarray*}
as claimed. As before, the upper bound can be derived in the same manner. 
{\hfill$\Box$}\bigskip 

\begin{lemma}\label{W11} We have $\varphi\in W^{1, 1}(]0, 1[)$. 
\end{lemma} 
{\bf Proof\,:} From Lemma \ref{le5}(a) 
we deduce that $\varphi'\in L^1(]0, 1-\Gamma[)$ and 
\[ \varphi(s)=\varphi(s_0)+\int_{s_0}^s\varphi'(\sigma)\,d\sigma \] 
holds for all $s, s_0\in ]0, 1-\Gamma[$. Hence $\varphi$ is absolutely 
continuous on $[0, 1-\Gamma]$. Due to Lemma \ref{le5}(b) 
the same conclusion can be drawn on $[1-\Gamma, 1]$. 
Thus the assertion follows from the following property of 
absolutely continuous functions: Assume that $\phi: [a, b]\to\R$ 
is continuous and, for some $c\in ]a, b[$, both restrictions 
of $\phi$ to $[a, c]$ and $[c, b]$ are absolutely continuous. 
Then $\phi$ is absolutely continuous on $[a, b]$. 
{\hfill$\Box$}\bigskip 

\begin{cor}\label{rothst} We have the bound 
\[ \|\varphi'\|_{L^1(]0, 1[)}\le C, \] 
where $C>0$ is independent of $\Gamma\in [1-2^{-12}, 1[$. 
\end{cor} 
{\bf Proof\,:} The bounds 
from Lemma \ref{le5}(a) and (b) yield 
\begin{eqnarray*} 
   \varphi'(s) & \le & \bigg(\frac{1}{\pi^2}\,\frac{1}{\Gamma\sqrt{s}}
   +\frac{1}{\pi^2}\frac{128}{(1-\Gamma)^{1/3}}\,\sqrt{s}\bigg)
   \,{\bf 1}_{\,[0, 1-\Gamma[}(s)
   \\ & & +\,\bigg(\frac{2+128\,(1-\Gamma)^{2/3}}{\pi^2\sqrt{s}}
   +\frac{1}{\Gamma^2\pi^2}\frac{1-\Gamma}{\sqrt{s-(1-\Gamma)}}\bigg)   
   \,{\bf 1}_{\,]1-\Gamma, 1]}(s).  
\end{eqnarray*}    
Therefore, 
\begin{eqnarray*}  
   \|\varphi'\|_{L^1(]0, 1[)} & \le & \int_0^{1-\Gamma}
   \bigg(\frac{1}{\pi^2}\,\frac{1}{\Gamma\sqrt{s}}
   +\frac{1}{\pi^2}\frac{128}{(1-\Gamma)^{1/3}}\,\sqrt{s}\bigg)\,ds 
   \\ & & +\,\int_{1-\Gamma}^1
   \bigg(\frac{2+128\,(1-\Gamma)^{2/3}}{\pi^2\sqrt{s}}
   +\frac{1}{\Gamma^2\pi^2}\frac{1-\Gamma}{\sqrt{s-(1-\Gamma)}}\bigg)\,ds
   \\ & \le & C_1\bigg(\sqrt{1-\Gamma}+(1-\Gamma)^{7/6}\bigg) 
   \\ & & +\,C_1
   \bigg([1+(1-\Gamma)^{2/3}](1-\sqrt{1-\Gamma})
   +(1-\Gamma)\sqrt{\Gamma}\bigg)
   \\ & \le & C_2,      
\end{eqnarray*} 
as claimed. 
{\hfill$\Box$}\bigskip 

\begin{cor}\label{koenstr} We have the pointwise bound 
\[ \bigg|\varphi'(s)-\frac{1}{\pi^2}\,\frac{1}{\sqrt{s}}\bigg|
   \le C\,(1-\Gamma)^{1/6}
   +C\,\frac{1-\Gamma}{\sqrt{s}}\,{\bf 1}_{[0, 1-\Gamma[}(s)
   +C\frac{1-\Gamma}{\sqrt{s-(1-\Gamma)}}\,{\bf 1}_{[1-\Gamma, 1[}(s) \] 
for $s\in [0, 1]\setminus\{1-\Gamma\}$, 
where $C>0$ is independent of $\Gamma\in [1-2^{-12}, 1[$. 
\end{cor} 
{\bf Proof\,:} If $s\in [0, 1-\Gamma[$, then (\ref{mudy}) 
and (\ref{2ndt}) show that 
\begin{eqnarray*} 
   \bigg|\varphi'(s)-\frac{1}{\pi^2}\,\frac{1}{\sqrt{s}}\bigg| 
   & \le & \frac{1}{\pi^2}\,\frac{1}{\sqrt{s}}\,
   \,\bigg|\frac{1}{\Gamma}-1\bigg|
   +\frac{1}{\pi^2}\,\bigg|\int_0^s\frac{\psi''(s-\sigma)}{\sqrt{\sigma}}\,d\sigma
   \bigg|
   \\ & \le & C_1\,\frac{1}{\sqrt{s}}\,(1-\Gamma)
   +C_1\,\frac{1}{(1-\Gamma)^{1/3}}\,\sqrt{s}
   \\ & \le & C_1\,\frac{1}{\sqrt{s}}\,(1-\Gamma)
   +C_1\,(1-\Gamma)^{1/6}. 
\end{eqnarray*} 
Similarly, if $s\in ]1-\Gamma, 1]$, we can use (\ref{2ablb}) and (\ref{flopl}) to get 
\begin{eqnarray*} 
   \bigg|\varphi'(s)-\frac{1}{\pi^2}\,\frac{1}{\sqrt{s}}\bigg| 
   & \le & \frac{1}{\pi^2}\,\frac{1}{\sqrt{s}}\,
   \,\bigg|\frac{1}{\Gamma}-1\bigg|
   +\frac{1}{\pi^2}\,\bigg|\int_{s-(1-\Gamma)}^s\frac{\psi''(s-\sigma)}{\sqrt{\sigma}}\,d\sigma
   \bigg|
   \\ & & +\,\frac{1}{\pi^2}\frac{1-\Gamma}{\Gamma^2\sqrt{s-(1-\Gamma)}}
   \\ & \le & C_2\,\frac{1}{\sqrt{s}}\,(1-\Gamma)
   +C_2\,(1-\Gamma)^{2/3}\,\frac{1}{\sqrt{s}}  
   +C_2\frac{1-\Gamma}{\sqrt{s-(1-\Gamma)}}
   \\ & \le & C_3\,(1-\Gamma)^{2/3}
   +C_3\frac{1-\Gamma}{\sqrt{s-(1-\Gamma)}}, 
\end{eqnarray*}
noting that $s\ge 1/2$ in this case. 
{\hfill$\Box$}\bigskip 

\begin{remark}\label{Gam1_2}  
{\rm The construction of $\varphi$ extends to the case $\Gamma=1$ as follows. 
Due to Remark \ref{Gam1_1} we have $\psi(s)=s$ for $s\in [0, 1]$. 
From (\ref{varphiform}) this implies 
\[ \varphi(s)=\frac{1}{\pi^2}
   \frac{d}{ds}\int_0^s\frac{s-\sigma}{\sqrt{\sigma}}\,d\sigma
   =\frac{4}{3\pi^2}
   \frac{d}{ds}\Big(s^{3/2}\Big)
   =\frac{2}{\pi^2}\,\sqrt{s} \]  
for $s\in [0, 1]$. 
}  
\end{remark} 

%%%%%%%%%%%%%%%%%%%%%%%%%%%%%%%%%%%%%%%%%%%%%%%%%%%%%%%%%%%%%%%%%%%%%%%%%%%%%%%%%%

\setcounter{equation}{0}
\section{The function $\Phi_\Gamma$ 
and end of the proof of Theorem \ref{ex_fGam}}
\label{PhiGam_sect} 

For $\Gamma\in [1-2^{-12}, 1[$ let $\varphi=\varphi_\Gamma$ 
be as in Lemma \ref{le4} and define $\Phi=\Phi_\Gamma: [-1, 0]\to\R$ by 
\[ \Phi(u)=\varphi(u+1),\quad u\in [-1, 0]. \] 
We will always extend $\Phi$ by zero outside of $[-1, 0]$. 
Writing out the various definitions, 
(\ref{zielgl}) says that 
\[ \int_0^s\frac{\varphi(\sigma)}{\sqrt{s-\sigma}}\,d\sigma
   -\int_0^{s-\chi_\Gamma(s)}\frac{\varphi(\sigma)}{\sqrt{s-\chi_\Gamma(s)-\sigma}}\,d\sigma
   =\frac{s}{\pi},
   \quad s\in [0, 1]. \] 
Shifting the second integral, this may be restated as 
\[ \int_0^s\frac{\varphi(\sigma)}{\sqrt{s-\sigma}}\,d\sigma
   -\int_{\chi_\Gamma(s)}^s\frac{\varphi(\sigma-\chi_\Gamma(s))}{\sqrt{s-\sigma}}\,d\sigma
   =\frac{s}{\pi},
   \quad s\in [0, 1]. \]
Thus if we put $r=\sqrt{1-s}$, we get 
\begin{equation}\label{rhoeq1}
   \frac{\pi}{1-r^2}\bigg[\int_0^{1-r^2}\frac{\varphi(y)}{\sqrt{(1-r^2)-y}}\,dy
   -\int_{\tilde{\chi}_\Gamma(r)}^{1-r^2}
   \frac{\varphi(y-\tilde{\chi}_\Gamma(r))}
   {\sqrt{(1-r^2)-y}}\,dy\bigg]=1,
   \quad r\in [0, 1[;
\end{equation}  
note that 
\begin{equation}\label{asbis} 
   \tilde{\chi}_\Gamma(r)=\chi_\Gamma(s)
   =\left\{\begin{array}{c@{\quad:\quad}c} 
   \Gamma\,\frac{1-r^2}{r^2} & r\in [\sqrt{\Gamma}, 1]
   \\[1ex] 1-r^2 & r\in [0, \sqrt{\Gamma}]\end{array}\right.
\end{equation} 
in terms of the variable $r$. Splitting up (\ref{rhoeq1}) 
into two cases according to (\ref{asbis}), we see that in fact 
(\ref{rhoeq1}) is equivalent to (\ref{rhoeq2}). In other words, 
having solved (\ref{rhoeq1}) for $\varphi$, we have shown that 
$\Phi$ is a solution to (\ref{rhoeq2}). Thus in order to finish the proof 
of Theorem \ref{ex_fGam}, we need to check that $\varphi$ has all the properties 
that have been listed at the beginning of Section \ref{reform_sect}. 
But those properties have been derived in Lemmas \ref{le4}, \ref{le5} 
and \ref{W11}. 
{\hfill$\Box$}\bigskip 

\noindent 
{\bf Proof of Corollary \ref{sing_cor}:} Once again we drop the index $\Gamma$, 
since $\Gamma$ will be fixed. First we consider the $5$-manifolds 
${\cal M}$ and ${\cal N}$ from (\ref{calM_def}) and (\ref{calN_def}), 
respectively. In the variables $(r, p_r, u)$ they are parametrized by 
\[ \{(r, p_r, u)\in [0, 1]\times\R\times [-1, 0]: u=-1, 
   -\gamma-p_r^2\le u\} \] 
and 
\[ \{(r, p_r, u)\in [0, 1]\times\R\times [-1, 0]: u=-\Gamma, 
   -\gamma-p_r^2\le u\}, \] 
and the action of ${\rm SO}(3)$, as described in Lemma \ref{intlem}(a) 
in the appendix; also cf.~(\ref{pic2d}) and recall (\ref{gamma_def}). 
According to (\ref{fGam}) and (\ref{konsu}) we have 
\[ f(x, v)=\frac{3}{4\pi}\,\varphi'(u+1)\,
   {\bf 1}_{\{\ell^2\le\Gamma\}},
   \quad (x, v)\in\R^3\times\R^3, \] 
and $\varphi'(s)\to\infty$ as $s\searrow 0$ or $s\searrow 1-\Gamma$ 
by Lemma \ref{le5}(a) and (b). This immediately implies 
(\ref{cor_c2}). Next, from the identity (\ref{mudy}) 
we deduce that the limit 
\[ l_\Gamma=\lim_{s\nearrow 1-\Gamma}\varphi'(s) \] 
exists and moreover that $l_\Gamma\to\infty$ as $\Gamma\nearrow 1$. 
As a consequence, we obtain (\ref{cor_c1}). The conclusion 
on the behavior in ${\cal D}\setminus ({\cal M}\cup {\cal N})$ 
can even be improved in view of Lemma \ref{le5}: in fact we have 
$f_\Gamma\in C^\infty({\cal D}\setminus ({\cal M}\cup {\cal N}))$. 
{\hfill$\Box$}\bigskip 

%%%%%%%%%%%%%%%%%%%%%%%%%%%%%%%%%%%%%%%%%%%%%%%%%%%%%%%%%%%%%%%%%%%%%%%%%%%%%%%%%%

\setcounter{equation}{0}
\section{The limit $\Gamma\to 1$: Proof of Theorem \ref{Gamto1}}

From (\ref{fGam}) and (\ref{konsu}) we recall that 
\begin{equation}\label{fGamord} 
   f_\Gamma(x, v)=\frac{3}{4\pi}\,\Phi'_\Gamma(\ell^2(x, v)-|x|^2-|v|^2)
   \,{\bf 1}_{\{\ell^2(x, v)\le\Gamma\}},
   \quad (x, v)\in\R^3\times\R^3,
\end{equation}  
where now we make the dependencies on $\Gamma$ explicit. 
First we will argue that $f_1=f_{{\rm Kurth}}$, 
where $f_1$ is $f_\Gamma$ for $\Gamma=1$. 
According to Remark \ref{Gam1_2}, 
we have $\varphi_\Gamma(s)=\frac{2}{\pi^2}\,\sqrt{s}$ for $s\in [0, 1]$. 
From the definition, this yields $\Phi_1(u)=\varphi_1(u+1)
=\frac{2}{\pi^2}\,\sqrt{u+1}$ for $u\in [-1, 0]$ and $\Phi_1(u)=0$ 
for $u\in\R\setminus [-1, 0]$, and hence in particular 
\begin{equation}\label{flei}
   \Phi'_1(u)=\frac{1}{\pi^2}\,\frac{1}{\sqrt{u+1}},\quad u\in ]-1, 0[. 
\end{equation}    
Comparing this to (\ref{fKdef}), we see that indeed $f_1=f_{{\rm Kurth}}$. 

For the $L^1$-estimate, let 
\[ \phi(x, v)={\rm sgn}\,(f_\Gamma(x, v)-f_{{\rm Kurth}}(x, v)), \] 
where we take ${\rm sgn}(0)=0$. Then $\phi$ is defined a.e.~and $\phi\in L^\infty(\R^3\times\R^3)$. 
Using (\ref{freg}) for $\Gamma$ and $\Gamma=1$, we get
\begin{eqnarray*} 
   \lefteqn{\int_{\R^3}\int_{\R^3} |f_\Gamma(x, v)-f_{{\rm Kurth}}(x, v)|\,dx\,dv} 
   \\ & = & \int_{\R^3}\int_{\R^3} (f_\Gamma(x, v)-f_{{\rm Kurth}}(x, v))\,\phi(x, v)\,dx\,dv
   \\ & = & \frac{3}{8\pi}\int_0^{\sqrt{\Gamma}} dr\,\frac{r^2}{1-r^2}
   \int_{|p_r|\le\sqrt{1-r^2}} dp_r\int_{-1}^{-r^2-p_r^2} du 
   \,\Phi'_\Gamma(u)\,I_\phi(r, p_r, \ell(r, p_r, u))
   \\ & & +\,\frac{3}{8\pi}\int_{\sqrt{\Gamma}}^1 dr\,\frac{r^2}{1-r^2}
   \int_{|p_r|\le\sqrt{1-\gamma}} dp_r\int_{-\gamma-p_r^2}^{-r^2-p_r^2} du
   \,\Phi'_\Gamma(u)\,I_\phi(r, p_r, \ell(r, p_r, u))
   \\ & & +\,\frac{3}{8\pi}\int_{\sqrt{\Gamma}}^1 dr\,\frac{r^2}{1-r^2}
   \int_{\sqrt{1-\gamma}\le |p_r|\le\sqrt{1-r^2}} dp_r  
   \int_{-1}^{-r^2-p_r^2} du\,\Phi'_\Gamma(u)\,I_\phi(r, p_r, \ell(r, p_r, u))
   \\ & & -\,\frac{3}{8\pi}\int_0^1 dr\,\frac{r^2}{1-r^2}
   \int_{|p_r|\le\sqrt{1-r^2}} dp_r\int_{-1}^{-r^2-p_r^2} du 
   \,\Phi'_1(u)\,I_\phi(r, p_r, \ell(r, p_r, u)). 
\end{eqnarray*} 
Since $|\phi(x, v)|\le 1$ a.e., we get $|I_\phi(r, p_r, \ell)|\le 8\pi^2$ 
by (\ref{Iphi_bd}). From Corollary \ref{rothst} it follows that 
\begin{eqnarray}\label{kureis} 
   \lefteqn{\int_{\R^3}\int_{\R^3} |f_\Gamma(x, v)-f_{{\rm Kurth}}(x, v)|\,dx\,dv} 
   \nonumber   
   \\ & \le & C_1\,\bigg|\int_0^{\sqrt{\Gamma}} dr\,\frac{r^2}{1-r^2}
   \int_{|p_r|\le\sqrt{1-r^2}} dp_r\int_{-1}^{-r^2-p_r^2} du 
   \,\Phi'_\Gamma(u)\,I_\phi(r, p_r, \ell(r, p_r, u))
   \nonumber   
   \\ & & \hspace{2em} -\,\int_0^1 dr\,\frac{r^2}{1-r^2}
   \int_{|p_r|\le\sqrt{1-r^2}} dp_r\int_{-1}^{-r^2-p_r^2} du 
   \,\Phi'_1(u)\,I_\phi(r, p_r, \ell(r, p_r, u))\,\bigg|
   \nonumber   
   \\ & & +\,C_1\int_{\sqrt{\Gamma}}^1 dr\,\frac{r^2}{1-r^2}
   \int_{|p_r|\le\sqrt{1-\gamma}} dp_r  
   \nonumber   
   \\ & & +\,C_1\int_{\sqrt{\Gamma}}^1 dr\,\frac{r^2}{1-r^2}
   \int_{\sqrt{1-\gamma}\le |p_r|\le\sqrt{1-r^2}} dp_r  
   \nonumber   
   \\ & \le & C_2\,\int_0^{\sqrt{\Gamma}} dr\,\frac{r^2}{1-r^2}
   \int_{|p_r|\le\sqrt{1-r^2}} dp_r\int_{-1}^{-r^2-p_r^2} du 
   \,|\Phi'_\Gamma(u)-\Phi'_1(u)|
   \nonumber   
   \\ & & +\,C_2\,\int_{\sqrt{\Gamma}}^1 dr\,\frac{r^2}{1-r^2}
   \int_{|p_r|\le\sqrt{1-r^2}} dp_r\int_{-1}^{-r^2-p_r^2} du\,\Phi'_1(u)
   \nonumber   
   \\ & & +\,C_2\int_{\sqrt{\Gamma}}^1 dr\,\frac{r^2}{\sqrt{1-r^2}}
   \nonumber   
   \\ & \le & C_3\,\int_0^{\sqrt{\Gamma}} dr\,\frac{r^2}{1-r^2}
   \int_{|p_r|\le\sqrt{1-r^2}} dp_r\int_{-1}^{-r^2-p_r^2} du 
   \,|\Phi'_\Gamma(u)-\Phi'_1(u)|
   \nonumber   
   \\ & & +\,C_3\int_{\sqrt{\Gamma}}^1\frac{dr}{\sqrt{1-r}}\,,
\end{eqnarray} 
where here and below the $C_i>0$ are independent of $\Gamma\in [1-2^{-12}, 1[$. 

Next observe that, by definition and by Corollary \ref{koenstr},  
\begin{eqnarray*} 
   |\Phi'_\Gamma(u)-\Phi'_1(u)|
   & = & \bigg|\varphi'_\Gamma(u+1)-\frac{1}{\pi^2}\,\frac{1}{\sqrt{u+1}}\bigg|
   \\ & \le & C_4\,(1-\Gamma)^{1/6}
   +C_4\,\frac{1-\Gamma}{\sqrt{u+1}}\,{\bf 1}_{[0, 1-\Gamma[}(u+1)
   \\ & & \,+\,C_4\frac{1-\Gamma}{\sqrt{u+1-(1-\Gamma)}}\,{\bf 1}_{[1-\Gamma, 1[}(u+1)
\end{eqnarray*}   
for a.e.~$u\in [-1, 0]$. In particular, we get  
\[ \int_{-1}^0 |\Phi'_\Gamma(u)-\Phi'_1(u)|\,du\le C_5(1-\Gamma)^{1/6}+C_5(1-\Gamma)
   \le C_6(1-\Gamma)^{1/6}. \] 
Going back to (\ref{kureis}), this yields 
\begin{eqnarray*} 
   \lefteqn{\int_{\R^3}\int_{\R^3} |f_\Gamma(x, v)-f_{{\rm Kurth}}(x, v)|\,dx\,dv} 
   \\ & \le & C_7(1-\Gamma)^{1/6}\,\int_0^{\sqrt{\Gamma}} dr\,\frac{r^2}{1-r^2}
   \int_{|p_r|\le\sqrt{1-r^2}} dp_r
   +C_7\sqrt{1-\sqrt{\Gamma}}
   \\ & \le & C_8(1-\Gamma)^{1/6}+C_8(1-\Gamma)^{1/2}
   \\ & \le & C_9(1-\Gamma)^{1/6}. 
\end{eqnarray*} 
This completes the proof of Theorem \ref{Gamto1}. 
{\hfill$\Box$}\bigskip 

%%%%%%%%%%%%%%%%%%%%%%%%%%%%%%%%%%%%%%%%%%%%%%%%%%%%%%%%%%%%%%%%%%%%%%%%%%%%%%%%%%

\setcounter{equation}{0}
\section{Periodic solutions $f_{\Gamma, \eps}$ close to $f_\Gamma$: 
Proof of Theorem \ref{per_gam}}
\label{Liling} 

We start with some auxiliary results. 

\begin{lemma}\label{thaif} 
Assume that ${\cal T}: I\times\R^3\times\R^3\to\R^3\times\R^3$ 
is a $C^1$-function such that the map ${\cal T}_t(x, v)={\cal T}(t, x, v)$ 
is a measure-preserving diffeomorphism of $\R^3\times\R^3$ 
for every $t\in I$. Then, given $f_\ast\in L^1(\R^3\times\R^3)$, 
the map 
\[ f: I\to L^1(\R^3\times\R^3), 
   \quad f(t)(x, v)=f_\ast({\cal T}_t(x, v)),
   \quad t\in I,\quad (x, v)\in\R^3\times\R^3, \] 
is continuous.   
\end{lemma} 
{\bf Proof\,:} Since ${\cal T}_t$ is measure-preserving 
for $t\in I$ and as $f_\ast\in L^1(\R^3\times\R^3)$, 
$f: I\to L^1(\R^3\times\R^3)$ is well-defined. 

To prove the continuity of $f$, we first consider the case 
where $f_\ast$ is continuous and has compact support $K\subset\R^3\times\R^3$. 
Let $(t_j)\subset I$ be such that $t_j\to t_0\in I$ as $j\to\infty$. 
Then we can fix an interval $J\subset\joinrel\subset I$ 
so that $t_j\in J$ as $j\in\N$. For $t\in I$ we denote by 
${\cal S}_t=({\cal T}_t)^{-1}$ the diffeomorphism of $\R^3\times\R^3$ 
that is inverse to ${\cal T}_t$. From the implicit function theorem 
we obtain that ${\cal S}$ is a $C^1$-function, where ${\cal S}(t, x, v)
={\cal S}_t(x, v)$. Then the set $\tilde{K}={\cal S}(J\times K)$ 
is compact in $\R^3\times\R^3$ and contains the support 
of $f(t)$ for every $t\in J$. Now if $(x, v)\in\R^3\times\R^3$, then 
\[ f(t_j)(x, v)=f_\ast({\cal T}(t_j, x, v))
   \to f_\ast({\cal T}(t_0, x, v))=f(t_0)(x, v) \]  
as $j\to\infty$ by the continuity of $f_\ast$ and ${\cal T}$. 
In addition, we note that 
\[ |f(t_j)(x, v)-f(t_0)(x, v)|\le 2\,{\|f_\ast\|}_\infty 
   {\bf 1}_{\tilde{K}}(x, v),
   \quad j\in\N,\quad (x, v)\in\R^3\times\R^3. \] 
Hence the dominated convergence theorem implies that 
$f(t_j)\to f(t_0)$ in $L^1(\R^3\times\R^3)$ as $j\to\infty$. 

Once we know the continuity of $f$ in this special case, 
we can proceed to the proof in the general case. 
If $f_\ast\in L^1(\R^3\times\R^3)$, then given $\eps>0$ 
there is a continuous function $f_{\ast, \eps}\in L^1(\R^3\times\R^3)$ 
of compact support such that ${\|f_\ast-f_{\ast, \eps}\|}_{L^1(\R^3\times\R^3)}
\le\eps/3$. Since ${\cal T}_t$ is measure-preserving, we have 
\[ {\|f(t)-f_\eps(t)\|}_{L^1(\R^3\times\R^3)}\le\frac{\eps}{3},
   \quad t\in I, \] 
where $f_\eps(t)=f_{\ast, \eps}\circ {\cal T}_t$. 
Let $t_0\in I$. Due to the special case there is $\delta>0$ such that 
$t\in I$ and $|t-t_0|\le\delta$ implies that 
\[ {\|f_\eps(t)-f_\eps(t_0)\|}_{L^1(\R^3\times\R^3)}\le\frac{\eps}{3}. \] 
As a consequence, we obtain 
\[ {\|f(t)-f(t_0)\|}_{L^1(\R^3\times\R^3)}\le\eps \] 
for $t\in I$ so that $|t-t_0|\le\delta$. 
{\hfill$\Box$}\bigskip 

\begin{proposition}\label{per_gen} Assume that $f_\ast=f_\ast(x, v)$ 
is a non-negative weak static solution such that 
\[ \rho_\ast=\rho_{{\rm Kurth}}, \] 
cf.~(\ref{hano}), where 
\[ \rho_\ast(x)=\int_{\R^3} f_\ast(x, v)\,dv,
   \quad {\rm a.e.}\,\,x\in\R^3. \] 
Suppose also that the essential support of $f_\ast$ 
is contained in $\{(x, v)\in\R^3\times\R^3: |x|\le 1, |v|\le 1\}$. 
In addition, let $\phi=\phi(t)$ for $t\in I$ 
be a solution of (\ref{Kep}), 
i.e.~$\ddot{\phi}=-\frac{1}{\phi^2}+\frac{1}{\phi^3}$ and $\phi>0$. Then 
\begin{equation}\label{fphi} 
   f_\phi(t, x, v)=f_\ast\Big(\frac{x}{\phi(t)}, 
   \phi(t)v-\dot{\phi}(t)x\Big),
   \quad t\in I,\quad (x, v)\in\R^3\times\R^3,
\end{equation}  
defines a weak solution. 
\end{proposition} 
{\bf Proof\,:} Let  
\[ {\cal T}: I\times\R^3\times\R^3\to\R^3\times\R^3,
   \quad {\cal T}(t, x, v)=\Big(\frac{x}{\phi(t)}, 
   \phi(t)v-\dot{\phi}(t)x\Big),
   \quad t\in I,\quad (x, v)\in\R^3\times\R^3. \] 
Then ${\cal T}$ is a $C^1$-function 
and every ${\cal T}_t={\cal T}(t, \cdot, \cdot)$ 
is a measure-preserving diffeomorphism of $\R^3\times\R^3$; 
note that in this case the inverses are given by 
\[ {\cal S}_t(y, w)=({\cal T}_t)^{-1}(y, w)
   =\Big(\phi(t)y, \frac{w}{\phi(t)}+\dot{\phi}(t)y\Big),
   \quad t\in I,\quad (y, w)\in\R^3\times\R^3. \] 
Hence we deduce from Lemma \ref{thaif} that $f: I\to L^1_+(\R^3\times\R^3)$ 
is continuous. Let $J\subset\joinrel\subset I$. 
Then $\min_{t\in J}\phi(t)>0$, $\max_{t\in J}\phi(t)<\infty$ 
and $\max_{t\in J}|\dot{\phi}(t)|<\infty$. Since we are assuming 
that the essential support of $f_\ast$ 
is contained in $\{(x, v)\in\R^3\times\R^3: |x|\le 1, |v|\le 1\}$, 
it follows that one can find a compact set $K_J\subset\R^3\times\R^3$ so that
$f(t, x, v)=0$ for all $t\in J$ and 
a.e.~$(x, v)\in (\R^3\times\R^3)\setminus K_J$. 
It remains to verify that (\ref{distr_eq}) is satisfied. 
Since $\rho_\ast=\rho_{{\rm Kurth}}$ by hypothesis, 
we observe that the definition (\ref{fphi}) of $f_\phi$ yields 
\[ \rho_{f_\phi}(t, x)
   =\int_{\R^3} f_\phi(t, x, v)\,dv
   =\frac{1}{\phi(t)^3}\,\rho_\ast
   \Big(\frac{x}{\phi(t)}\Big)
   =\frac{1}{\phi(t)^3}\,\rho_{{\rm Kurth}}
   \Big(\frac{x}{\phi(t)}\Big), \]     
and hence
\begin{equation}\label{hostr} 
   U_{f_\phi}(t, x)
   =\frac{1}{\phi(t)}\,U_{{\rm Kurth}}\Big(\frac{x}{\phi(t)}\Big);
\end{equation}   
cf.~(\ref{rhoeps}) and (\ref{Ueps}). 
As the change of variables $(t, {\cal T}_t(x, v))=(t, y, w)$ 
is measure-preserving, we calculate 
for $\varphi\in {\cal D}(I\times\R^3\times\R^3)$ 
\begin{eqnarray}\label{J1}  
   B & := & \int_I\int_{\R^3}\int_{\R^3} f_\phi(t, x, v)
   \,\Big(\partial_t\varphi(t, x, v)+v\cdot\nabla_x\varphi(t, x, v)
   -\nabla_x U_{f_\phi}(t, x)\cdot\nabla_v\varphi(t, x, v)\Big)
   \,dt\,dx\,dv
   \nonumber   
   \\ & = & \int_{\R^3}\int_{\R^3} dy\,dw
   \,f_\ast(y, w)\int_I dt\, 
   \Big[\partial_t\varphi(t, {\cal S}_t(y, w))
   +\Big(\frac{w}{\phi(t)}+\dot{\phi}(t)y\Big)
   \cdot\nabla_x\varphi(t, {\cal S}_t(y, w))
   \nonumber   
   \\ & & \hspace{12.5em} -\,\nabla_x U_{f_\phi}(t, \phi(t)y)
   \cdot\nabla_v\varphi(t, {\cal S}_t(y, w))\Big]. 
\end{eqnarray} 
The essential support of $f_\ast(y, w)$ 
is contained in $\{(y, w)\in\R^3\times\R^3: |y|\le 1, |w|\le 1\}$. 
If $|y|<1$, then $|\phi(t)y|<\phi(t)$ 
and therefore $\nabla_x U_{f_\phi}(t, \phi(t)y)
=\frac{y}{\phi(t)^2}$ in view of (\ref{hostr}) 
and (\ref{UQkurth}). Accordingly, (\ref{J1}) reduces to   
\begin{eqnarray}\label{J2}  
   B & = & \int_{\R^3}\int_{\R^3} dy\,dw
   \,f_\ast(y, w)\int_I dt\, 
   \Big[\partial_t\varphi(t, {\cal S}_t(y, w))
   +\Big(\frac{w}{\phi(t)}+\dot{\phi}(t)y\Big)
   \cdot\nabla_x\varphi(t, {\cal S}_t(y, w))
   \nonumber   
   \\ & & \hspace{12.5em} -\,\frac{y}{\phi(t)^2}
   \cdot\nabla_v\varphi(t, {\cal S}_t(y, w))\Big].    
\end{eqnarray} 
Next we observe that 
\begin{eqnarray*} 
   \frac{\partial}{\partial t}\Big[\varphi(t, {\cal S}_t(y, w))\Big]
   & = & \partial_t\varphi(t, {\cal S}_t(y, w))
   +\dot{\phi}(t)y\cdot\nabla_x\varphi(t, {\cal S}_t(y, w))
   \\ & & +\,\Big(\ddot{\phi}(t)y
   -\frac{\dot{\phi}(t)}{\phi(t)^2}w\Big)
   \cdot\nabla_v\varphi(t, {\cal S}_t(y, w)).
\end{eqnarray*}  
Hence we obtain from (\ref{J2}) the relation 
\begin{eqnarray*}  
   B & = & \int_{\R^3}\int_{\R^3} dy\,dw\,f_\ast(y, w)
   \\ & & 
   \times\,\int_I dt\,\Big[
   \frac{\partial}{\partial t}\Big[\varphi(t, {\cal S}_t(y, w))\Big]
   +\frac{w}{\phi(t)}
   \cdot\nabla_x\varphi(t, {\cal S}_t(y, w))
   \\ & & \hspace{5em}
   +\,\Big(-\ddot{\phi}(t)y-\frac{y}{\phi(t)^2}
   +\frac{\dot{\phi}(t)}{\phi(t)^2}w\Big)
   \cdot\nabla_v\varphi(t, {\cal S}_t(y, w))\Big]
   \\ & = & \int_I\frac{dt}{\phi(t)^2}
   \int_{\R^3}\int_{\R^3} dy\,dw\,f_\ast(y, w)\,\Big[
   \phi(t)w\cdot\nabla_x\varphi(t, {\cal S}_t(y, w))
   \\ & & \hspace{14em} -\,\Big(\frac{y}{\phi(t)}
   -\dot{\phi}(t) w\Big)
   \cdot\nabla_v\varphi(t, {\cal S}_t(y, w))\Big],  
\end{eqnarray*} 
where we used that $\varphi$ vanishes at the endpoints 
of the time interval along with (\ref{Kep}). 
Next we introduce, at fixed $t\in I$, 
the function $\tilde{\varphi}\in {\cal D}(\R^3\times\R^3)$ by means of 
\[ \tilde{\varphi}(y, w)=\varphi(t, {\cal S}_t(y, w)),
   \quad(y, w)\in\R^3\times\R^3. \] 
Since $f_\ast$ is a weak static solution, we know by definition that 
\[ \tilde{B}(t):=\int_{\R^3}\int_{\R^3} dy\,dw\,f_\ast(y, w)
   \,\Big[w\cdot\nabla_y\tilde{\varphi}(y, w)
   -y\cdot\nabla_w\tilde{\varphi}(y, w)\Big]=0. \]
As 
\[ B=\int_I\frac{dt}{\phi(t)^2}\,\tilde{B}(t), \] 
it follows that $B=0$, which in turn shows 
that (\ref{distr_eq}) is verified. 
{\hfill$\Box$}\bigskip 

\noindent 
{\bf Proof of Theorem \ref{per_gam}\,:} Since every $f_\Gamma$ 
is a non-negative weak static solution by Theorem \ref{ex_fGam}, 
we can simply apply Proposition \ref{per_gen}. 
{\hfill$\Box$}\bigskip 

\begin{remark}\label{Kep_rem}   
{\rm The solutions of (\ref{Kep}) are all global in time 
and it is possible to derive an almost explicit formula for them. 
In fact, $\phi$ is the norm of a solution of Kepler's problem 
with angular momentum $1$. More precisely, $\phi(t)=|x(t)|$, where 
\[ \ddot{x}=-\frac{x}{|x|^3},\quad x\in\R^2\setminus\{0\},
   \quad |x\wedge\dot{x}|=1. \] 
In the elliptic case, i.e., for negative energy 
$H(x, \dot{x})=\frac{1}{2}\,|\dot{x}|^2-\frac{1}{|x|}=-\frac{1}{2(1-e^2)}$, 
we have 
\[ x(t)=\frac{1}{1-e^2}(\cos u-e, \sin u), \]  
where $e\in [0, 1[$ and the relation among $t$ and $u$ is 
\[ u-e\sin u=\frac{1}{(1-e^2)^{3/2}}\,(t-t_0). \] 
Then we get 
\[ \phi(t)=|x(t)|=\frac{1}{1-e^2}\sqrt{1+e^2-2e\cos u} \] 
and hence obtain periodic solutions of period $2\pi(1-e^2)^{3/2}$. 
As a further remark, the connection to the solutions $\phi_\eps$ 
of (\ref{Kep}) from (\ref{phieps_IVP}) is given by $1-e^2=\frac{1}{1-\eps^2}$. 

In the parabolic and hyperbolic cases, it is possible to derive 
analogous formulas in terms of hyperbolic functions. In this way 
one can generate additional (non-periodic) solutions of the Vlasov-Poisson system 
for all time. 
}  
\end{remark}

%%%%%%%%%%%%%%%%%%%%%%%%%%%%%%%%%%%%%%%%%%%%%%%%%%%%%%%%%%%%%%%%%%%%%%%%%%%%%%%%%%

\setcounter{equation}{0}
\section{An infinite-dimensional continuum of equilibria}
\label{infinite_sect} 

In this section the following general principle will be applied. 

\begin{lemma}\label{conv_lem} 
Suppose that $f_1, \ldots, f_N$ are weak solutions on the time interval $I\subset\R$ 
and for every $J\subset\joinrel\subset I$ there exists a compact set $K_J\subset\R^3\times\R^3$ 
such that the essential support of every ${f_i|}_J$ is contained in $K_J$. 
In addition, let the $f_i$ have a common density, i.e., 
\[ \rho_{f_i}=\rho_\ast,\quad i=1, \ldots, N, \] 
for some function $\rho_\ast: I\to L^1(\R^3)$. Then any convex combination 
\[ f=\sum_{i=1}^N\lambda_i f_i,\quad\lambda_i\in [0, 1],
   \quad\sum_{i=1}^N\lambda_i=1, \] 
is also a weak solution on $I$. 
\end{lemma} 
{\bf Proof\,:} Since 
\[ \rho_f=\rho_{\sum_{i=1}^N\lambda_i f_i}
   =\sum_{i=1}^N\lambda_i\rho_{f_i}=\rho_\ast, \] 
we get $U_f=U_{f_i}$ for $i=1, \ldots, N$. Hence it remains to multiply (\ref{distr_eq}) 
for $f_i$ by $\lambda_i$ and add $\sum_{i=1}^N$ to see that (\ref{distr_eq}) also holds for $f$. 
{\hfill$\Box$}\bigskip 

\noindent 
{\bf Proof of Theorem \ref{infi_thm}\,:} First we assert that, given 
\[ 1-10^{-12}\le\Gamma_1<\ldots<\Gamma_N\le 1, \] 
the corresponding weak static solutions $f_{\Gamma_i}$, $i=1, \ldots, N$, 
as constructed in Theorem \ref{ex_fGam} are linearly independent. To establish this claim, 
suppose that 
\begin{equation}\label{afmtg} 
   \sum_{i=1}^N\alpha_i f_{\Gamma_i}=0
\end{equation}  
for some coefficients $\alpha_i\in\R$. 
We have ${\cal N}_{\Gamma_i}\cap {\cal N}_{\Gamma_j}=\emptyset$ for $i\neq j$, 
cf.~(\ref{calN_def}), and furthermore all ${\cal M}_{\Gamma_i}, {\cal N}_{\Gamma_i}\subset {\cal D}_{\Gamma_i}
\subset\{(x, v)\in\R^3\times\R^3: |x|\le 1, |v|\le 1\}$ are compact.  
Since (\ref{cor_c2}) from Corollary \ref{sing_cor} shows that 
$f_{\Gamma_i}(x, v)\to\infty$ as $(x, v)\nearrow{\cal N}_{\Gamma_i}$, 
(\ref{afmtg}) is only possible, if $\alpha_i=0$ for $i=1, \ldots, N$. 

As all $f_{\Gamma_i}$ have their essential support contained 
in $K=\{(x, v)\in\R^3\times\R^3: |x|\le 1, |v|\le 1\}$, 
Lemma \ref{conv_lem} implies that 
\[ {\cal E}(\Gamma_1, \ldots, \Gamma_N)
   =\Big\{\sum_{i=1}^N\lambda_i f_{\Gamma_i}: 
   \lambda_i\in [0, 1],\quad\sum_{i=1}^N\lambda_i=1\Big\} \] 
is a set of weak static solutions of dimension $N-1$. We can also consider 
the very large set 
\[ {\cal E}=\bigcup_{N=1}^\infty\,\,\bigcup_{1-10^{-12}\le\Gamma_1<\ldots<\Gamma_N\le 1}
   {\cal E}(\Gamma_1, \ldots, \Gamma_N) \] 
of weak static solutions, all of whose essential supports are in $K$. 
By Definition \ref{weak_sol_def}, it follows that then also ${\cal C}=\overline{{\cal E}}$ 
is a set of weak solutions with their essential supports contained in $K$, 
where the closure is taken in $L^1(\R^3\times\R^3)$. 
In addition, we have $\rho_{f_\ast}=\rho_{{\rm Kurth}}$ 
for every $f_\ast\in {\cal C}$.  
Lastly, applying the construction of periodic solutions (Proposition \ref{per_gen}), 
it follows that every weak static solution $f_\ast\in {\cal C}$ 
is surrounded by time-periodic weak solutions. 
{\hfill$\Box$}\bigskip 

%%%%%%%%%%%%%%%%%%%%%%%%%%%%%%%%%%%%%%%%%%%%%%%%%%%%%%%%%%%%%%%%%%%%%%%%%%%%%%%%%%

\setcounter{equation}{0}
\section{Appendix}
\label{append_sect} 

\begin{lemma}\label{intlem} 
(a) We consider the compact Lie group ${\rm SO}(3)$ and its 
Haar measure $dR$ such that $\int_{{\rm SO}(3)} dR=8\pi^2$. 
If $g\in L^1(\R^3\times\R^3)$, then 
\begin{equation}\label{intform} 
   \int_{\R^3}\int_{\R^3} g(x, v)\,dx\,dv
   =\int_0^\infty dr\int_{\R} dp_r\int_0^\infty d\ell\,\ell\,I_g(r, p_r, \ell),
\end{equation}  
where 
\[ I_g(r, p_r, \ell)=\int_{{\rm SO}(3)} 
   g\left(R\left(\begin{array}{c} 0 \\ 0 \\ r\end{array}\right),
   R\left(\begin{array}{c} \ell/r \\ 0 \\ p_r\end{array}\right)\right)\,dR. \]   
If $g$ is spherically symmetric and we identify 
$g(x, v)=\tilde{g}(|x|, x\cdot v/|x|, |x\wedge v|)$, 
then $I_{gh}=\tilde{g} I_h$ for general $h=h(x, v)\in L^\infty(\R^3\times\R^3)$.  
In particular, if $g$ is spherical symmetric, then $I_g=8\pi^2\tilde{g}$. 
\smallskip 

\noindent 
(b) If $r>0$ and $g(re_3, \cdot)\in L^1(\R^3)$, then 
\[ \int_{\R^3} g(re_3, v)\,dv
   =\frac{1}{2r^2}\int_{\R} dp_r\int_0^\infty d\ell^2\int_0^{2\pi} d\vartheta 
   \,g(re_3, v(r, p_r, \ell, \vartheta)), \] 
where $v(r, p_r, \ell, \vartheta)=(\frac{\ell}{r}\cos\vartheta, \frac{\ell}{r}\sin\vartheta, p_r)$. 
\end{lemma} 
{\bf Proof\,:} (a) For any $(x, v)\in\R^3\times\R^3$ such that $x\wedge v\neq 0$ 
there is a unique $R\in {\rm SO}(3)$ such that 
\[ Re_3=\frac{x}{|x|}=\frac{x}{r},\quad Re_2=-\frac{x\wedge v}{|x\wedge v|}=-\frac{x\wedge v}{\ell}. \] 
Then, with $p_r=\frac{x\cdot v}{r}$, 
\[ rp_r x-r^2 v=(x\cdot v)x-|x|^2 v=x\wedge (x\wedge v)
   =-(r Re_3)\wedge (\ell Re_2)=-r\ell\,R(e_3\wedge e_2)=-r\ell\,Re_1, \] 
and hence 
\[ v=p_r\frac{x}{r}+\frac{\ell}{r}\,Re_1=p_r Re_3+\frac{\ell}{r}\,Re_1
   =R\left(\begin{array}{c} \ell/r \\ 0 \\ p_r\end{array}\right), \] 
in addition to 
\[ x=R\left(\begin{array}{c} 0 \\ 0 \\ r\end{array}\right). \] 
The map 
\[ [0, \infty[\times\R\times [0, \infty[\times {\rm SO}(3) 
   \ni (r, p_r, \ell, R)\mapsto (x, v)
   =\left(R\left(\begin{array}{c} 0 \\ 0 \\ r\end{array}\right), 
   R\left(\begin{array}{c} \ell/r \\ 0 \\ p_r\end{array}\right)\right)
   \in\R^3\times\R^3 \] 
parametrizes $\R^3\times\R^3$, up to a set of measure zero.  
Using Euler angles to parametrize ${\rm SO}(3)$ and to write out $dR$, 
it can be shown that $dx\,dv=dr\,dp_r\,d\ell\,\ell\,dR$, 
which yields (\ref{intform}). For the last assertion, we observe 
\[ g\left(R\left(\begin{array}{c} 0 \\ 0 \\ r\end{array}\right),
   R\left(\begin{array}{c} \ell/r \\ 0 \\ p_r\end{array}\right)\right)
   =\tilde{g}(r, p_r, \ell), \] 
since 
\[ \left|R\left(\begin{array}{c} 0 \\ 0 \\ r\end{array}\right)\right|=r,
   \quad R\left(\begin{array}{c} 0 \\ 0 \\ r\end{array}\right)
   \cdot R\left(\begin{array}{c} \ell/r \\ 0 \\ p_r\end{array}\right)=rp_r,
   \quad R\left(\begin{array}{c} 0 \\ 0 \\ r\end{array}\right)
   \wedge R\left(\begin{array}{c} \ell/r \\ 0 \\ p_r\end{array}\right)
   =R\left(\begin{array}{c} 0 \\ \ell \\ 0\end{array}\right). \] 
(b) This follows by direct integration. 
{\hfill$\Box$}\bigskip 

%%%%%%%%%%%%%%%%%%%%%%%%%%%%%%%%%%%%%%%%%%%%%%%%%%%%%%%%%%%%%%%%%%%%%%%%%%%%%%%%%%%%

\end{document}